\documentclass[5p,twocolumn,10pt,times]{elsarticle}
\usepackage{amsmath}
\usepackage{bm}
\usepackage{times}
\usepackage{amsmath}
\usepackage{amssymb}
\usepackage{mathptmx}
\usepackage[stretch=10]{microtype}
\usepackage{comment}
\usepackage[hidelinks]{hyperref}
\usepackage{listings,xcolor}
\usepackage{adjustbox}
\usepackage[framed,numbered]{matlab-prettifier}
\usepackage{xcolor}
\usepackage[caption=false,font=normalsize,labelfont=sf,textfont=sf]{subfig}
\usepackage{paralist}

\usepackage[color=yellow]{todonotes}
\definecolor{lightGreen}{RGB}{173,223,138}
\definecolor{lightBlue}{RGB}{166,206,227}
\definecolor{lightYellow}{RGB}{255,238,170}

\def \bezier {B{\'e}zier}
\def \spectral   {\textsc{Spectral}}
\def \spectralPE {\textsc{Spectral\! PE}}
\def\COMMENT #1 {}

\begin{document}
\baselineskip11pt
\lstset{language=Matlab}
\begin{frontmatter}

\title{Spectral mesh-free quadrature for planar regions bounded by rational parametric curves}

\author[1,2]{David Gunderman\corref{cor1}}
\cortext[cor1]{Corresponding author}
\ead{david.gunderman@colorado.edu}
\author[2]{Kenneth Weiss}
\author[3]{John A. Evans}
\address[1]{Department of Applied Mathematics, University of Colorado Boulder, Boulder, CO 80309, USA}
\address[2]{Lawrence Livermore National Laboratory, 7000 East Avenue, Livermore, CA 94550, USA}
\address[3]{Ann and H.J. Smead Department of Aerospace Engineering Sciences, University of Colorado Boulder, Boulder, CO 80309, USA}

\begin{abstract} 
This work presents spectral, mesh-free, Green's theorem-based numerical quadrature schemes for integrating functions over planar regions bounded by rational parametric curves. 
Our algorithm proceeds in two steps:
\begin{inparaenum}[(1)]
  \item We first find intermediate quadrature rules for line integrals along the region's boundary curves corresponding to Green's theorem.
  \item We then use a high-order quadrature rule to compute the numerical antiderivative of the integrand along a coordinate axis, which is used to evaluate the Green's theorem line integral.
\end{inparaenum} 
We present two methods to compute the intermediate quadrature rule. The first is spectrally accurate (it converges faster than any algebraic order with respect to number of quadrature points) and is relatively easy to implement, but has no guarantee of polynomial exactness.  The second guarantees exactness for polynomial integrands up to a pre-specified degree $k$ with an \textit{a priori}-known number of quadrature points and retains the convergence properties of the first, but is slightly more complicated. The quadrature schemes have applications to computation of geometric moments, immersogeometric analysis, conservative field transfer between high-order meshes, and initialization of multi-material simulations with rational geometry. We compare the quadrature schemes produced using our method to other methods in the literature and show that they are much more efficient both in terms of number of quadrature points and computational time.  We provide an open-source implementation of the algorithm in MATLAB. 
\end{abstract}

\begin{keyword} 
    Quadrature, Rational, NURBS, Spectral Convergence
\end{keyword}
\end{frontmatter}


\section{Introduction}

Quadrature over arbitrary regions bounded by rational parametric curves is important in a variety of applications including computer-aided design (CAGD) and computer-aided manufacturing (CAM), in which geometric objects are typically represented in terms of their boundaries (BREPs) in a basis of non-uniform rational B-splines (NURBS). 
Efficient and accurate computation of integrated quantities, such as geometric moments and mass, are often key ingredients in the design process.

\begin{figure}
  \centering
	\vspace{-.7cm}
  \includegraphics[width=.85\linewidth]{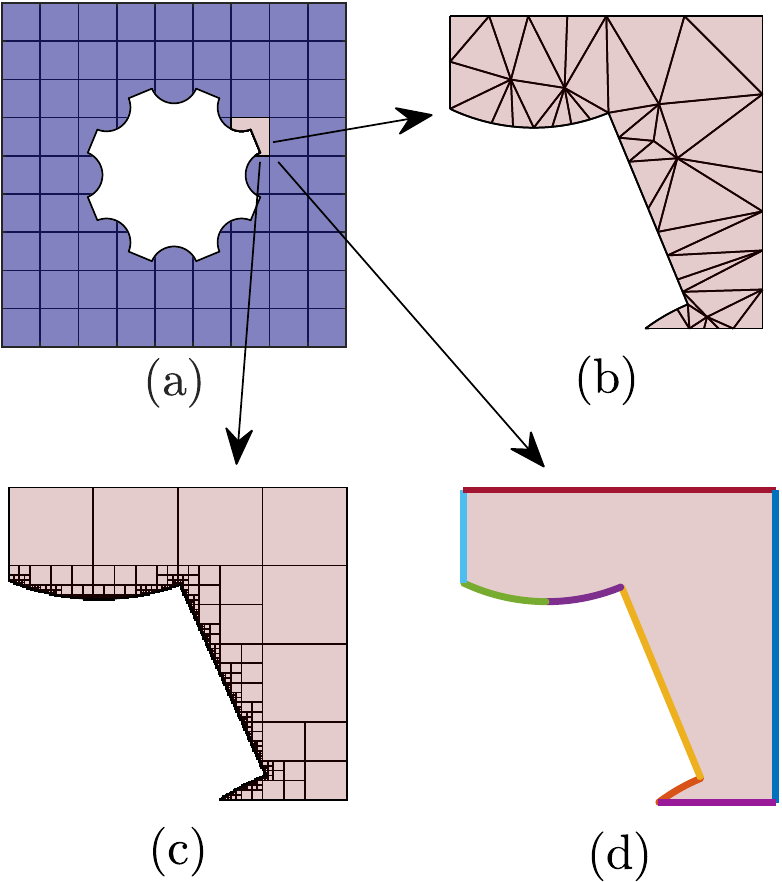}
  \caption{
    \textup{(a)} In immersogeometric analysis, background cells intersecting the immersed boundary require special quadrature rules. 
    \textup{(b)} Rational meshing takes the geometry into account, but requires expensive pre-processing. 
    \textup{(c)} Adaptive quadtree integration converges at most linearly. 
    \textup{(d)} Our boundary-based method takes the geometry into account without expensive meshing. }
  \vspace{-1cm}
  \label{fig:immersogeometric_example}
\end{figure}

\begin{figure*}
\includegraphics[width=\linewidth]{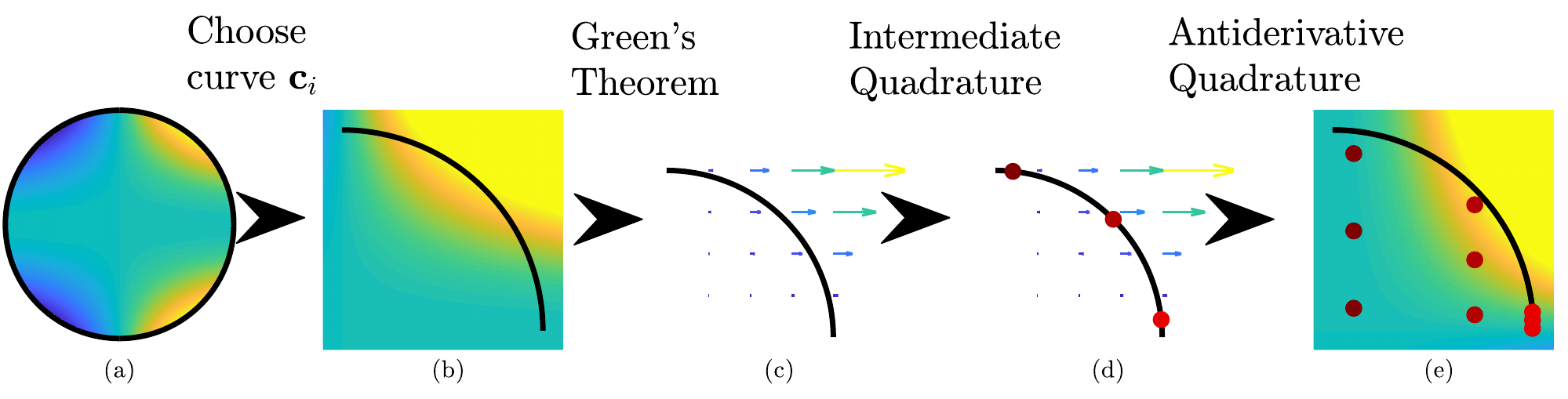}
\caption{Overview of our quadrature algorithms for integration over a region bounded by rational curves. 
         The integrand function $f(x,y)$ is given as a color plot in (a), (b), and (e),
	     while its antiderivative $A_f(x,y)$ is given as a colored vector plot in (c) and (d). 
	     The intermediate quadrature points are notated by the shaded red dots in (d). 
	     The final quadrature points, produced by evaluating the antiderivative function using antiderivative quadrature at the red dots in (d), 
	     are shown as correspondingly shaded dots in (e).}
      
\label{fig:overview}
\end{figure*}

While CAD has relied on NURBS as a basis for decades, rational parameterization of problem geometry has also recently become more popular in computer-aided analysis in an effort to more closely couple the design/analysis pipeline. In isogeometric analysis~\cite{hughes2005isogeometric, bazilevs2010isogeometric, borden2011isogeometric}, NURBS parameterizations are used to represent the interior of the objects. However, as only boundaries are typically represented in CAD, this requires the computation of a high-quality interior parameterization in the form of a mesh. Thus, recently, there has been a push to combine isogeometric analysis with immersed boundary treatment using immersogeometric analysis \cite{kamensky2015immersogeometric, schillinger2012isogeometric}.  In immersogeometric analysis, the basis functions are defined with respect to a background grid which is immersed within the geometry of interest.  This avoids the need to construct a high-quality interior parameterization, but requires integration over the intersections of the background grid with the geometry, as shown in Figure~\ref{fig:immersogeometric_example}.

There already exist relatively efficient integration strategies for traditional high-order parametric elements (e.g.\ triangles/tetrahedra or quads/hexahedra) which occur in isogeometric analysis using parametric mappings of high-order simplicial and tensor product-type quadrature rules~\cite{hughes2010efficient}. However, these quadrature rules cannot be seamlessly or efficiently applied to arbitrary rational regions due to the inherent difficulties in high-order meshing of arbitrary geometries. In contexts such as immersogeometric analysis~\cite{krishnamurthy2010accurate}, as well as in cutFEM \cite{burman2010fictitious}, field transfer between high-order meshes~\cite{anderson2018high}, initialization of volume of fluid fields for rational shapes~\cite{verschaeve2011high, dyadechko2005moment}, and moment-fitting~\cite{mousavi2010generalized}, efficient integration over arbitrary rational regions (such as intersections of conic sections) is fundamental. 

The need for efficient quadrature schemes for arbitrary planar rational regions could potentially be met by the rich body of literature on efficient integration over arbitrary planar regions. However, when considering rational geometry, these methods either require expensive preprocessing in the form of meshing or only achieve low-order convergence due to their low-order approximation of the geometry. Our algorithm is able to achieve spectral convergence (i.e. faster than any algebraic order with respect to the number of quadrature points) without the need for expensive meshing. 

Existing methods for integrating over arbitrary 2D geometries fall into two general categories: (1) domain decomposition methods and (2) Green's theorem methods. Domain decomposition methods generally rely on two steps: tessellation using a set of basic elements (such as a triangulation or a decomposition into nested rectangles) and integration over these basic elements. Another recently introduced method is the class of correction-based methods, in which geometry is approximated by linear boundaries, then a geometric correction term is applied \cite{thiagarajan2014adaptively, scholz2017first}. 

Low-order domain decomposition methods tessellate the region using linear polygons (quadtrees, triangulation, and quadrangulation being the most common cases). The tessellation and integration steps are both very robust in the case of linear polygonal elements and have been studied extensively. The robustness is paid for in low-order convergence \cite{flusser1999calculation, wu2001new, singer1993general, yang1996fast}, which cannot be overcome by increasing the quadrature order, since the underlying geometric approximation is low-order.

Higher-order or exact domain decomposition methods are generally less robust and require more preprocessing time, because they attempt to tessellate the geometry using high-order or rational parametric elements to achieve higher-order convergence. The tessellation problem over planar regions becomes a high-order 
2D meshing problem, which has been studied extensively~\cite{roca2011defining, sherwin2002mesh, engvall2016isogeometric, engvall2017isogeometric, engvall2018mesh}. Once the parametric elements have been defined, integrals over the elements are computed using well-known Gaussian quadrature rules for high-order parametric mappings of triangles or rectangles. High-order (and even spectral, in the case of exact meshing) convergence can be attained, because both the quadrature scheme and the geometric approximation can be made high-order.

The other common method for quadrature over complicated geometries avoids meshing the interior by using Green's theorem to transform double integrals over the domain into a sum of line integrals over its boundaries. These methods have been employed extensively to increase computational efficiency of integration in a variety of situations, including for linear polygonal regions, regions defined by polynomial parametric boundaries, and implicitly-defined regions~\cite{chin2015numerical, sudhakar2013quadrature,saye2015high, olshanskii2016numerical}.

A key limitation of Green's theorem-based methods is that they require some parametric description of the boundary. However, an approximation of the boundary using splines or other high-order parametric curves can be used to guarantee high-order convergence~\cite{sommariva2006meshless,sommariva2009gauss,santin2011algebraic}, while still avoiding the potential difficulties of meshing. If boundaries are given parametrically (as is the case for NURBS objects), this extra information is available without any extra approximation, and should be fully exploited to attain the most efficient quadrature rules possible -- a fact which we take advantage of in our scheme.
Another popular belief about Green's theorem methods is that they can only be used on integrands for which analytic antiderivatives are available, and these methods have been used in this case for decades~\cite{sheynin2003moment, li1993moment}. However, recently it has become more common to use Green's theorem-based methods for arbitrary integrands~\cite{saye2015high,sommariva2009gauss,jonsson2017cut}, since numerical antiderivatives can be efficiently computed to high precision using a high-order quadrature rule, such as Gaussian quadrature -- another fact which we take advantage of in our scheme.
%

In this paper, we consider the problem of integrating arbitrary functions over planar regions bounded by rational curves, with special emphasis on those that may be difficult to efficiently mesh with simple parametric elements. 
While specialized integration strategies based on Green's theorem have been previously published in the literature for regions bounded by linear segments, polynomial curves, and implicit curves, to the best of our knowledge, no such optimized integration method has been proposed for regions bounded by rational curves. 


More formally, we develop two numerical algorithms for computing quadrature rules for integrals of the form 
\begin{equation*}
  \iint_{\Omega} f(x,y) dA,
\end{equation*}
over an arbitrary planar region $\Omega$ bounded by a set of $n_c$ rational curves of degrees $\{m_i\}_{i=1}^{n_c}$, 
\begin{equation*}
\partial{\Omega} = \displaystyle\cup_{i=1}^{n_c} \mathbf{c}_i:
\end{equation*}
 \begin{equation*}
  \mathbf{c}_i(s) = 
    \begin{cases} x_i(s) \\
                  y_i(s)
    \end{cases} 
    1 \leq i \leq n_c,
  \end{equation*} 
\noindent
for $0\leq s \leq 1$.  


The first algorithm produces quadrature rules which, for smooth integrands, converge to the correct integral exponentially quickly with respect to the number of quadrature points. Exponential convergence is a particular class of spectral convergence (in which error decreases faster than any algebraic order), so we refer to this algorithm as \spectral\ quadrature. The second algorithm produces quadrature rules which retain the fast convergence properties of the first, but are also exact (up to machine precision) for polynomial integrands up to a pre-specified degree $k$, so we refer to it as \spectralPE\ (spectral, polynomially exact) quadrature.

As outlined in Figure~\ref{fig:overview}, our algorithms both consist of the same two primary steps, which are applied to each component curve $\mathbf{c}_i$ individually:
\begin{inparaenum}[(1)]
  \item we find appropriate intermediate quadrature rules for line integrals corresponding to Green's theorem, taking into account the degrees of the corresponding boundary curves and
  \item we use a high-order quadrature rule to compute the numerical antiderivative of the integrand along a coordinate axis, which is used to evaluate the Green's theorem line integral.
\end{inparaenum}
In the \spectral\ method, we use Gaussian quadrature for both quadrature rules, the order of which can be increased to converge faster than any algebraic order with respect to the total number of quadrature points. In the \spectralPE\ method, we use a rational quadrature rule instead of Gaussian quadrature for the intermediate quadrature rule. These methods are described in detail in Section~\ref{sec:algorithm}.

In addition to their spectral convergence, our algorithms are also computationally efficient.
The cost for producing the \spectral\ quadrature rules scales linearly with the number of points in each of the underlying Gaussian quadrature rules. The cost of producing a \spectralPE\ quadrature rule that is exact for polynomials up to degree $k$ scales as $\mathcal{O}(n_c m^3 k^2)$, where $m$ is the highest degree of any of the $n_c$ curves, $m=\max_{i\leq n_c}(m_i)$. Typically, the degrees, $m_i$, of the component curves are all equal. We derive the number of quadrature points and the computational costs in Section~\ref{sec:algorithm}.

For smooth integrands over arbitrary rational regions, we show on numerical examples in Section~\ref{sec:results} that both the \spectral\ and \spectralPE\ methods are much more efficient than existing methods both in terms of number of quadrature points and computational time. We compare our two methods to five methods used in the literature: quadtree integration~\cite{krishnamurthy2010accurate}, linear meshing with high-order Gaussian integration~\cite{engwirda2014locally}, rational meshing with high-order Gaussian integration~\cite{engvall2016isogeometric}, linear approximation of the boundary combined with a Green's theorem approach~\cite{sommariva2009gauss}, and a cubic spline approximation of the boundary combined with a Green's theorem approach~\cite{sommariva2009gauss}. 



\paragraph{Contributions}
The main contributions of our work include:
\begin{compactitem}
    \item We present algorithms which compute the locations and weights of \spectral\ quadrature rules for integration of arbitrary smooth integrands and \spectralPE\ (spectral, polynomially exact) quadrature rules for exact (up to machine precision) integration of polynomials up to a given degree $k$ over a region bounded by rational curves. Both algorithms converge faster than any algebraic order, as we demonstrate on a variety of test shapes and integrands.
    \item We derive a formula for the number of quadrature points needed in the \spectralPE\ quadrature rule for exact (up to machine precision) integration of polynomials up to degree $k$ over an arbitrary region bounded by rational curves. This is analogous to Gaussian quadrature for polynomial regions, for which an exact number of quadrature points can be calculated to integrate polynomials up to a pre-specified degree in 1D.
    \item A claimed difficulty of our type of approach in the literature is that it requires antiderivatives. We show that these are available numerically and easy to calculate exactly for a broad class of useful integrands, including any polynomial. Moreover, we can compute highly accurate numerical antiderivatives of an even broader class of integrands, as can be seen in our experiments in Figure~\ref{fig:nonpolynomial_results}.
    \item We compare our proposed quadrature schemes to five other commonly-used quadrature schemes in the literature, which represent the state of the art for integrating over rational regions, and show that the proposed algorithms produce quadrature rules which are significantly more efficient both in terms of the number of quadrature points and the computation time needed to achieve a given integration accuracy.
\end{compactitem}

\paragraph{Outline}
The remainder of the paper is organized as follows: 
In Section~\ref{sec:algorithm}, we introduce the new algorithms for the \spectral\ and \spectralPE\ quadrature rules for planar regions bounded by rational curves. 
We present numerical results on a variety regions and for a variety of integrands with comparison methods in Section~\ref{sec:results} 
and conclude in Section~\ref{sec:discussion} with a discussion of our results and ideas for future work.

\section{Numerical Integration over Rational Regions}
\label{sec:algorithm}

We first consider a genus-zero planar region $\Omega$ with a closed, connected, oriented boundary loop $\Gamma$ such that the normal to $\Gamma$ points toward the interior of the region $\Omega$. 
We take the boundary loop $\Gamma$ to be composed of $n_c$ rational oriented component curves $\{\mathbf{c}_i\}_{i=1}^{n_c}$, where the $\mathbf{c}_{i}$ are given in the rational Bernstein-\bezier\ basis:

\begin{equation*}
  \mathbf{c}_i(s) = 
    \begin{cases} x_i(s) = \frac{\displaystyle\sum_{j=0}^{m_i} w_{i,j} x_{i,j} B^{m_i}_j(s)}{\displaystyle\sum_{j=0}^{m_i} w_{i,j} B^{m_i}_j(s)} \\
                  y_i(s) =  \frac{\displaystyle\sum_{j=0}^{m_i} w_{i,j} y_{i,j} B^{m_i}_j(s)}{\displaystyle\sum_{j=0}^{m_i} w_{i,j} B^{m_i}_j(s)}
    \end{cases} 
    1 \leq i \leq n_c,
\end{equation*}
   for $0\leq s \leq 1$ such that $\Gamma = \cup_{i=1}^{n_c} \mathbf{c}_i$, where $m_i$ are the degrees of the component curves, $w_{i,j}$ are the control weights, and $x_{i,j},y_{i,j}$ are the control points.  The Bernstein-\bezier\ basis is defined as 
  \begin{equation}
  B_j^m(s) = {m\choose{j}} (1-s)^{j}s^{m-j},
  \end{equation} and has a multitude of useful properties~\cite{farin2002curves,farouki2012}. For rational curves, the control weight matrix $w_{i,j}$ has at least one entry such that $w_{i,j} \neq 1$ for some $i,j$. We denote the maximum degree of the component curves as $m= \max_{i \leq n_c} (m_i)$.
\begin{figure}
  \centering
  \includegraphics[width=.8\linewidth]{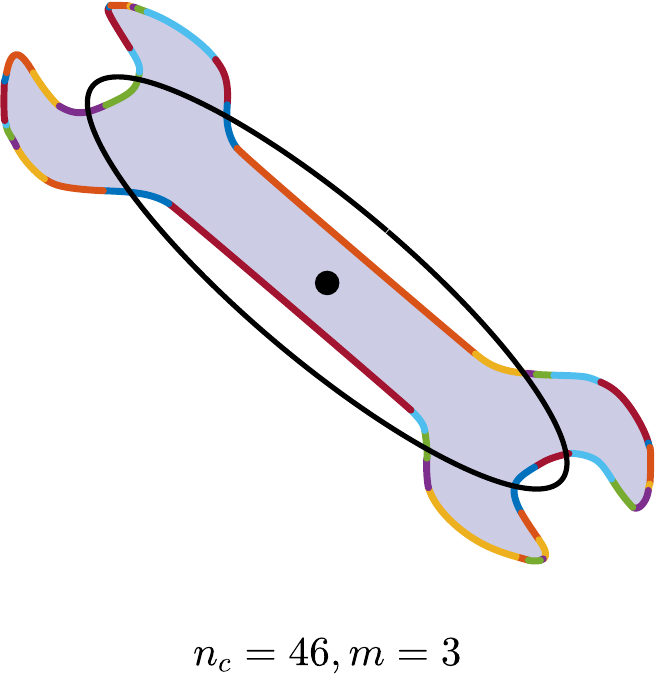}
  \caption{
    A wrench-shaped region bounded by rational curves over which polynomials can be integrated to machine precision using the \spectralPE\ algorithm. 
    This example illustrates the object's first few geometric moments, as computed by the \spectralPE\ algorithm:
    The grey area represents its area (0$^{th}$ moment), 
    the black dot represents its centroid (1$^{st}$ moment), 
    and the ellipse represents its moment of inertia (2$^{nd}$ moment). }
  \label{fig:example_region}
\end{figure}
We denote the bounding box of the control points of the region $\Omega$ as
\begin{equation*}
\tilde{\Omega} = [\min_{i,j}\{x_{i,j}\}, \max_{i,j}\{x_{i,j}\} ] \times [\min_{i,j}\{y_{i,j}\}, \max_{i,j}\{y_{i,j}\} ].
\end{equation*} 

Both the \spectral\ and \spectralPE\ algorithms compute quadrature points $\{(x_l,y_l)\}_{l=1}^{n_q}$ and weights $\{w_l\}_{l=1}^{n_q}$ such that the error
\begin{equation*}
E = \left|\iint _{\Omega} f(x,y) dx dy - \sum_{l=1}^{n_q} w_l f(x_l,y_l)\right| < \mathcal{O}(e^{-cn_q}),
\end{equation*}  
decreases exponentially with respect to number of quadrature points for some constant $c$ for smooth integrands in the bounding box of $\Omega$, $f(x,y) \in \mathcal{C}^\infty(\tilde{\Omega})$,
where $n_q$ is the number of quadrature points. This is because they are based on one-dimensional high-order quadrature rules, which can be shown to have exponential convergence (c.f. \cite{davis2007methods}).
If the function has only $k$ continuous derivatives,  $f(x,y) \in \mathcal{C}^k (\tilde{\Omega})$, then the scheme attains the maximum algebraic order of convergence possible, $E = \mathcal{O}(n_q^{-k+1})$.

Moreover, because it uses rational quadrature for the intermediate quadrature rule, the \spectralPE\ quadrature rule is exact up to machine precision, $\epsilon$, for polynomials $p(x,y)$ of maximal degree $k$ or less with an \textit{a priori} known number of quadrature points $n_q$, which, as we derive in Section~$\ref{sec:modification}$, is:
\begin{equation}
	n_q = \left\lceil{\frac{k+1}{2}}\right\rceil\sum_{i=1}^{n_c} \left(m_i(k + 3)+1\right).
\end{equation}

Figure~\ref{fig:example_region} shows an example of a wrench-shaped region composed of $n_c=46$~component curves, all of which are rational and of degree $m_i=3$.
The \spectralPE\ quadrature scheme requires $n_q = 1472$ quadrature points to exactly integrate polynomials up to degree $k=2$ and $n_q = 1748$ quadrature points to exactly integrate polynomials up to degree $k=3$. This works out to $32$ and $38$ quadrature points per component curve, respectively, for $k=2$ and $k=3$.
Note that the \spectralPE\ quadrature scheme still retains the spectral convergence property for arbitrary smooth integrands, so fewer quadrature points could be used if machine precision is not necessary.

It is important here to emphasize the limitation that the integrand $f(x,y)$ must be smoothly extendable to the bounding box $\tilde{\Omega}$ of the integration region $\Omega$. This is trivially possible when the integrand is polynomial, and is also possible in immersed finite element methods such as immersogeometric analysis and the finite cell method. In these methods, the integrand is defined on a background cell, and the bounding box of the integration region is contained within this background cell. An example is the pink cell in Figure~\ref{fig:immersogeometric_example}. In cases when the integrand cannot be smoothly extended, moment-fitting can be used with the \spectralPE\ algorithm to produce high-order quadrature rules with quadrature points within the domain of integration. We briefly describe moment fitting in~\ref{app:momentFitting}.

The general idea behind the present algorithms is to use Green's theorem to transform the integral over the interior region $\Omega$ into a sum of line integrals along each curve $\mathbf{c}_i$,

\begin{align}
\iint_\Omega f(x,y) dxdy &= \int_{\partial \Omega} A_f(x,y) dx \nonumber \\
&= \sum_{i=1}^{n_c} \int_{\mathbf{c}_i} A_f(x,y) dx, \label{eq:green}
\end{align}
where $A_f(x,y)$ is the $y$-antiderivative of the integrand.

The algorithms are overviewed graphically in Figure~\ref{fig:overview} and, for each component curve $\mathbf{c}_i$, can both be summarized in two basic steps:
\begin{enumerate}[(1)]
  \item Compute a set of intermediate quadrature points and weights for the Green's Theorem line integral corresponding to the particular curve $\mathbf{c}_i$, and
  \item Compute quadrature rules to evaluate the antiderivative function numerically at each intermediate quadrature point.
\end{enumerate}
The two quadrature rules' weights are then combined in a tensor-product fashion to produce the quadrature weights for each of the final quadrature points. The \spectral\ and \spectralPE\ algorithms only differ in  the implementation for the first step in the above two-step process. We describe the \spectral\ algorithm initially, then describe how the \spectralPE\ algorithm differs.

\subsection{Intermediate Quadrature Rule}
\label{sec:intermediate}
The first step in the algorithm is to convert the area integral inside the domain $\Omega$ into line integrals along each of the curves, $\{\mathbf{c}_{i}\}_{i=1}^{n_c}$, and to find appropriate intermediate quadrature rules for each of these line integrals. 
For the purposes of this section, we assume that the antiderivative function
\begin{equation}
A_f(x,y) = \int_{C}^y f(x,t)dt, \label{eq:antiderivative}
\end{equation} 
can be evaluated to high precision at any $(x,y)$ point on any of the component curves $\mathbf{c}_i$. 
We explain our approach for this in the next section (Section~\ref{sec:antiderivative}).

Because the curves are given parametrically, the physical space integrals appearing on the right-hand side of Equation~$\eqref{eq:green}$ can be transformed into the parametric space of each component curve:
\begin{equation}
\int_{\mathbf{c}_i} A_f(x,y) dx = \int_0^1 A_f(x_i(s),y_i(s)) \frac{d x_i(s)}{ds} ds. \label{eq:green_par}
\end{equation}
Each of these integrals can be evaluated using a high-order 1D quadrature scheme, such as Gaussian quadrature with $Q_i$-points, which yields a set of $Q_i$ quadrature points $\{s_{i,q}\}_{q=1}^{Q_i}$ and weights $\{\gamma_{i,q}\}_{q=1}^{Q_i}$. The intermediate quadrature points in $\mathbb{R}^2$ can then be calculated as $(x_{i,q},y_{i,q})= \left(x_i(s_{i,q}),y_i(s_{i,q})\right)$ for each $q=1,\ldots Q_i$. This evaluation can be performed stably for rational Bernstein-\bezier\ curves using de Casteljau's algorithm \cite{farin2002curves}.

The intermediate quadrature rule then reads as
\begin{equation}
\iint _{\Omega} f(x,y) dx dy \approx \sum_{i=1}^{n_c} \sum_{q=1}^{Q_i} \gamma_{i,q} A_f(x_{i,q},y_{i,q}) \frac{d x_{i}}{ds}(s_{i,q}). \label{eq:intermediate}
\end{equation}

The resulting integrand in each of the line integrals is composed of two parts, one which requires evaluation of the integrand's antiderivative and one which requires evaluation of the derivative of a \bezier\ curve. Derivatives of rational or polynomial \bezier\ curves can be efficiently and stably computed to machine precision using de Casteljau's algorithm. As long as the antiderivative is evaluated to high enough precision, the intermediate quadrature rule will converge exponentially fast as $Q_i$ increases \cite{davis2007methods}.  

Next, we describe how to evaluate the antiderivative.

\subsection{Antiderivative Quadrature Rule}
\label{sec:antiderivative}

A key ingredient in the evaluation of the interior function of the intermediate quadrature rule in Equation~$\eqref{eq:intermediate}$ is the evaluation of the antiderivative. A high-order quadrature rule, such as Gaussian quadrature with $P$ points, can be used to evaluate the integrand with spectral convergence as $P$ is increased~\cite{davis2007methods}. In our implementation, for simplicity, we set $P=Q_i$. Note that for some integrands, such as polynomial integrands, the antiderivative can alternatively be evaluated symbolically, which could garner computational efficiency benefits.

The antiderivative quadrature rule will then have quadrature points 
$\{x_{i,q}, y_{i,q,\zeta}\}_{\zeta=1}^{P}$, and weights $\{\gamma_{i,q,\zeta}\}_{\zeta=1}^{P}$. For a particular curve $\mathbf{c}_i$ and intermediate quadrature point $(x_{i,q},y_{i,q})$, the antiderivative can be evaluated as 
\begin{equation}
A_f(x_{i,q},y_{i,q}) = \int_{C}^{y_{i,q}} f(x_{i,q},y) dy \approx \sum_{\zeta=1}^{P} \gamma_{i,q,\zeta}\ f(x_{i,q},y_{i,q,\zeta}). \label{eq:num_antiderivative}
\end{equation}
Note that the $x_{i,q}$ do not depend on the index $\zeta$, since the antiderivative is only with respect to $y$. If the function $f(x,y)$ is a polynomial of degree $k$, then Equation~$\eqref{eq:num_antiderivative}$ will be exact (up to machine precision) with $P= \lceil{\frac{k+1}{2}\rceil}$ points. If the function is smooth, then convergence to the correct integral will be spectral with respect to the number of quadrature points as $P$ is increased.

Although the arbitrary integration constant $C$ in the calculation of $A_f(x,y)$ has no bearing on the validity of Equation~$\eqref{eq:green}$, it can affect the locations and weights of the resulting quadrature scheme as well as the floating point stability, if it chosen to be very far away from the integration domain. We choose $C$ according to the strategy described in ~Section~\ref{sec:green_implement}. 

\subsection{Full Quadrature Rule}
Finally, the intermediate quadrature rule and the numerical antiderivative quadrature rule can be combined to obtain a spectral full quadrature scheme:
\begin{equation}
\iint _{\Omega} f(x,y) dx dy \approx \sum_{i=1}^{n_c} \sum_{q=1}^{Q_i} \sum_{\zeta=1}^{P}  w_{i,q,\zeta} f(x_{i,q},y_{i,q,\zeta}) ,
\end{equation}
where $w_{i,q,\zeta} = \gamma_{i,q} \gamma_{i,q,\zeta} \frac{d x_{i}}{ds}(s_{i,q})$. The derivative terms $\frac{dx_{i}}{ds}$ can be computed efficiently using de Casteljau's algorithm and can be viewed as geometric correction terms to the component Gaussian quadrature weights. The total number of quadrature points used in the final quadrature rule will scale with the number of quadrature points used in each component quadrature rule, $n_q = P\sum_{i=1}^{n_c} Q_i$. If we use $P=Q_i$ and all $Q_i$ equal, then the number of quadrature points will be $n_q = n_c P^2$.

\subsection{A Modified Algorithm for Polynomial Exactness} 
\label{sec:modification}
In some situations, such as calculation of geometric moments and moment-fitting techniques, it is of particular importance to efficiently calculate integrals of polynomial functions over arbitrary rational regions. Calculating these integrals exactly with an \textit{a priori} number of quadrature points can remove the need for potentially expensive adaptive quadrature that attempts to converge to the solution through successively increased number of quadrature points. With a modification to the intermediate quadrature step of the \spectral\ algorithm presented above, we can achieve exactness up to machine precision for polynomial integrands over these regions. We call this quadrature scheme \spectralPE\ (spectral, polynomially exact). Exactness can be achieved by guaranteeing that the two components steps of intermediate quadrature and antiderivative quadrature are exact.

The antiderivative quadrature rule is already exact for polynomial integrands $f(x,y)$ up to a pre-specified degree $k$ if Gaussian quadrature is used with $P=\lceil{\frac{k+1}{2}\rceil}$ points. The intermediate quadrature rule, on the other hand, has no guarantees of exactness for polynomial integrands if Gaussian quadrature is used, since both the antiderivative function $A_f(x(s),y(s))$ and the parametric derivative term $\frac{dx}{ds}$ are rational, rather than polynomial, functions. However, a rational quadrature rule exact for an appropriate class of intermediate rational functions can be used for the intermediate quadrature rule in order to achieve polynomial exactness (up to machine precision). See \ref{app:detailed_example} for a worked-out example of the \spectralPE\ algorithm on a circular domain.
 
\subsubsection{Review of Rational Quadrature}
 Exact quadrature formulas for particular classes of functions such as polynomials have been studied for centuries. Quadrature rules which are exact for polynomials are generally formulated in terms of classes of orthogonal polynomials~\cite{gautschi1999orthogonal}. Two of the most commonly used quadrature rules are Gaussian quadrature and Fej\'er (i.e. Clenshaw-Curtis) quadrature, which are based on the orthogonal Legendre and Chebyshev polynomials, respectively. Exact quadrature formulas for classes of nonpolynomial functions have also been studied in the literature, and are typically based on generalizations of the classical orthogonal polynomials. The first study of orthogonal rational functions was made by Djirbashian in the 1960's, a review of which can be found in~\cite{djrbashian1990survey}. A connection to exact Gauss-type quadrature rules for rational functions with known real poles was first made by Gautschi~\cite{gautschi1993Gauss}. 

More recently, Bultheel, Deckers, and Van Deun have published a series of papers with very efficient algorithms for computing Fej\'er-type quadrature rules exact for rational functions with known poles outside the unit interval in linear time with respect to the total number of poles of the rational functions for which exactness is guaranteed~\cite{deckers2008rational, van2008algorithm, deckers2009computing, deckers2017algorithm}. The algorithm only requires knowing the locations and multiplicities of the poles of the rational functions of interest. 

Next, we describe an algorithm which computes locations and multiplicities of the poles of the rational functions which must be integrated exactly in order to achieve polynomial exactness.

\subsubsection{Poles of Intermediate Rational Functions}
\label{subsec:poles}
For a particular curve $\mathbf{c}_i$, the intermediate rational function which must be integrated exactly is the one that arises as the product of the antiderivative function $A_f(x_i(s))$ and the derivative term $\frac{dx(s)}{ds}$, as seen in Equation~$\eqref{eq:green_par}$. As we shall prove, the poles of the intermediate rational function $A_f(x_i(s)) \frac{dx(s)}{ds}$ do not in fact depend on the function $f(x,y)$ but only on its polynomial degree $k$. In fact, the poles of the intermediate rational function are simply the poles of the original curve $\mathbf{c}_i$ but with the multiplicity of each pole multiplied by $k+3$. To prove this, consider the two terms in the product separately. For clarity, we call $\{p_{i,j}\}_{j=1}^{m_i}$ the poles of the curve $\mathbf{c}_i$, defined as the zeros of the weight function 
\begin{equation*}
w_i(s) =\displaystyle\sum_{j=0}^{m_i} w_{i,j} B^{m_i}_j(s).
\end{equation*}

By the quotient rule, the derivative of a rational function has the effect of doubling the multiplicity of each of its poles. Therefore, $\frac{d x_i(s)}{ds}$ will have the same poles $\{p_{i,j}\}_{j=1}^{m_i}$ as $x_i(s)$, but each will have multiplicity of $2$. 

Since $x_i(s)$ and $y_i(s)$ are each rational functions with poles $\{p_{i,j}\}_{j=1}^{m_i}$ and $A_f(x,y)$ is a polynomial, $A_f(x(s),y(s))$ will also have the poles $\{p_{i,j}\}_{j=1}^{m_i}$. However, each pole's multiplicity will be $k+1$ times as large as the multiplicity it had for $x_i(s)$ or $y_i(s)$. Assuming that $f(x,y)$ is a degree $k$ polynomial and, expressing it in the monomial basis as $f(x,y) = \sum_{k_x+k_y\leq k} a_{k_x,k_y} x^{k_x} y^{k_y}$ for some coefficients $a_{k_x,k_y}$, then 
\begin{equation*}
\resizebox{.99\linewidth}{!}{\parbox{\linewidth}{
\begin{align*}
&A_f(x_i(s),y_i(s)) = \sum_{k_x+k_y\leq k} a'_{k_x,k_y} x_i(s)^{k_x} y_i(s)^{k_y+1} \\
&\hspace{-3mm}= \sum_{k_x+k_y \leq k} \left(a'_{k_x,k_y} \left(\frac{\sum_{j=0}^{m_i} w_{i,j}x_{i,j}B^{m_i}_j(s)}{\sum_{j=0}^{m_i} w_{i,j} B^{m_i}_j(s)}\right)^{k_x} \left(\frac{\sum_{j=0}^{m_i} w_{i,j}y_{i,j}s^j}{\sum_{j=0}^{m_i} w_{i,j} B^{m_i}_j(s)}\right)^{k_y+1} \right)\\
&\hspace{-3mm}= \sum_{k_x+k_y \leq k} \left(a'_{k_x,k_y}\frac{\left(\sum_{j=0}^{m_i} w_{i,j}x_{i,j}B^{m_i}_j(s)\right)^{k_x} \left(\sum_{j=0}^{m_i} w_{i,j}y_{i,j}B^{m_i}_j(s)\right)^{k_y+1}}{\left(\sum_{j=0}^{m_i} w_{i,j} B^{m_i}_j(s)\right)^{k_x+k_y+1}} \right),
\end{align*}
}}
\end{equation*}
where $a'_{k_x,k_y}=a_{k_x,k_y}/(k_y+1)$ and $k_x + k_y$ attains its maximum value of $k$, since $f(x,y)$ is a degree $k$ polynomial.

Combining the above statements, the integrand $A_f(x(s),y(s)) \frac{d x(s)}{d s}$ in Equation~$\eqref{eq:green}$ will have the same poles as the component curves did, namely $\{p_{i,j}\}_{j=1}^{m_i}$. However, each of these poles will have a multiplicity of $k+3$, where $k+1$ of them come from the antiderivative calculation and $2$ of them come from the derivative term. 

Therefore, in order to find a rational intermediate quadrature rule exact for rational functions with the required poles, we must compute the set of poles $\{p_{i,j}\}_{j=1}^{m_i}$ of each component curve $\{\mathbf{c}_i\}_{i=1}^{n_c}$. Because the rational curves are given in terms of control points and control weights, the poles of the rational curves $\{\mathbf{c}_i\}_{i=1}^{n_c}$ can be easily calculated as the roots of the corresponding control weight (i.e. denominator) polynomial, $w_i(s)=\sum_{j=0}^{m_i} w_{i,j} B^{m_i}_j(s)$.

There is a wide body of literature on stable univariate polynomial complex root finding. The roots of a polynomial in the monomial basis can be found by computing the eigenvalues of the polynomial's companion matrix, $A_i$ \cite{edelman1995polynomial}. By converting the polynomial $w_i(s)$ to the monomial basis $w_i(s) = \sum_{j=0}^{m_i} a_{i,j} s^j$ (see ~Section~\ref{sec:conversion}), the companion matrix can be written as 
\begin{align*}
A_i = \begin{bmatrix}
1 & 0 & 0 & \ldots & -\frac{a_{i,1}}{a_{i,m_i}}\\
1 & 1 & 0 & \ldots & -\frac{a_{i,2}}{a_{i,m_i}} \\
0 & 1 & 1 & \ldots & -\frac{a_{i,3}}{a_{i,m_i}} \\
\vdots & \vdots & \vdots & \ddots & \vdots \\
0 & 0 & 0 & \ldots & -1.
\end{bmatrix}
\end{align*}
The eigenvalues of $A_i$ are the roots of $w_i(s)$. Computing the eigenvalues of this $(m_i+1)\times (m_i+1)$ matrix is generally an $\mathcal{O}(m_i^3)$ operation. Finding the roots for all of the rational component curves therefore scales as $\mathcal{O}\left(n_c m^3\right)$ with the maximum degree of the component rational curves $m=\max_{i\leq n_c} (m_i)$, since finding the roots of the control weight polynomial $w_i(s) = \displaystyle\sum_{j=0}^{m_i} w_{i,j} B^{m_i}_j(s)$ requires $\mathcal{O}(m_i^3)$ operations. However, for relatively low degree rational curves, the cost of this root-finding step is minimal and the degree of the component rational curves is often taken to be constant in applications.

\subsubsection{Rational Intermediate Quadrature}
\label{sec:rat_intermediate}
Once the rational poles, $\{p_{i,j}\}_{i=1}^{m_i}$ for each $\mathbf{c}_i$ and their $k+3$ multiplicities are known, the remaining step is to calculate an exact quadrature rule for functions with these poles in the parametric space $0 \leq s \leq 1$ of each $\{\mathbf{c}_{i}\}_{i=1}^{n_c}$. Note that for any polynomial portions of the boundary, a simple Gaussian quadrature rule will be sufficient to achieve an exact intermediate quadrature rule. For the rational portions of the boundary, we must use an exact rational quadrature rule. We choose to use a linear-time algorithm developed in~\cite{deckers2017algorithm}. Implementation details are given in Section~\ref{sec:rational_implementation}. The quadrature points and weights from this algorithm can then be substituted for the Gaussian quadrature points and weights in Section~\ref{sec:intermediate}.

\begin{figure}
\begin{adjustbox}{width=.9\linewidth,keepaspectratio,margin=1.5em 0em 1.5em 0em}
\begin{lstlisting}[style=Matlab-editor,numbers=left]
f = @(x,y) getIntegrand();
kDeg=getDegree(f);
Bd=loadShape();
for i=1:getNumCurves(Bd)
  % Compute intermediate quadrature rule
  if use_SpectralPE
    poles=roots(getWeights(Bd(i)))  %2.5.1
    intPoles=repmat(poles,kDeg+3);  %2.4.2
    [intpts intwts]= 
         calcIntQRule(intPoles);    %2.5.3
  else
    [intpts intwts]= calcIntQRule();  %2.1
  end
  
  % Compute antiderivative quadrature rule
  for j=1:length(intpts)
    [Afx,Afy]=Bd{i}(intpts);
    [ADpts, ADwts]=                   %2.2
         calcAntiderivativeQRule([0,Afx]);
    finalptsx(j,:)=ADpts;
    len = length(finalptsx)
    finalptsy(j,:)=repmat(Afy,len,1); 
    finalwts(j,:)=intwts(j)*ADwts;
  end
  quadrules{i}=[finalptsx,
                finalptsy,
                finalwts ];           %2.3
end
\end{lstlisting}
\end{adjustbox}
\caption{Matlab pseudocode for the quadrature algorithms given in Section~\ref{sec:algorithm}. The function \lstinline{calcIntQRule} is the only difference between the \spectral\ algorithm and the \spectralPE\ algorithm. The first uses a typical high-order quadrature rule, such as Gaussian quadrature. The second uses a quadrature rule exact for rational functions with poles given by \lstinline{intPoles}. Line-by-line comments refer to sections of this article.}
%
\label{fig:alg_pseudocode}
\end{figure}

 The resulting exact quadrature rule will have  $n_q$ quadrature points:
 \begin{equation}
   n_q = \left\lceil{\frac{k+1}{2}}\right\rceil\sum_{i=1}^{n_c} \left(m_i(k + 3)+1\right),
 \end{equation}
where $n_c$ is the number of component curves, 
       $m_i$ are the degrees of the component curves $\mathbf{c}_i$, 
       and $k$ is the maximal degree of the polynomial integrand.
 Pseudocode for this algorithm is provided in Figure~\ref{fig:alg_pseudocode}.

\begin{figure}
\begin{minipage}{.235\textwidth}
\begin{center}
\includegraphics[width=\linewidth]{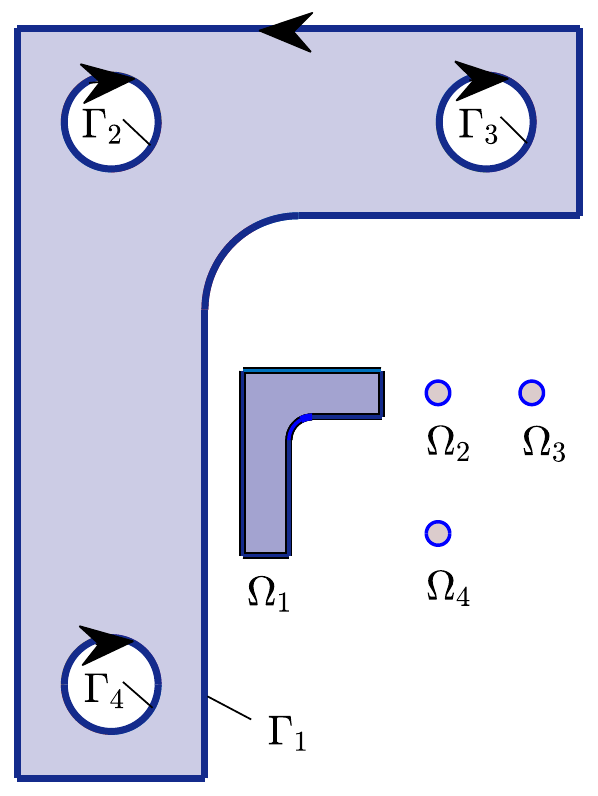}\\
\large{(a)}
\end{center}
\end{minipage}~
\begin{minipage}{.25\textwidth}
\begin{center}\includegraphics[width=\linewidth]{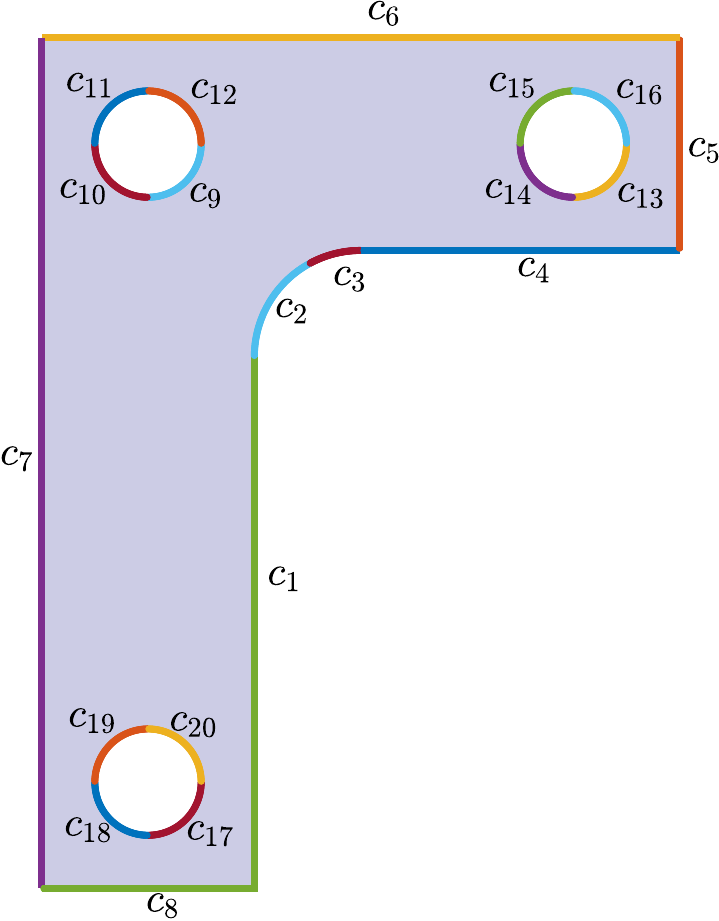}\\
\large{(b)}
\end{center}
\end{minipage}
\caption{An example of a non genus-zero region with $M=4$ boundary loops and $n_c=20$ component curves. In this case, the maximum degree of the component curves is $m=2$. In (a), the boundary curves $\{\Gamma_\alpha\}_{\alpha=1}^4$ and boundary regions $\{\Omega_\alpha\}_{\alpha=1}^4$ are labeled. In (b), the component curves $\{\mathbf{c}_i\}_{i=1}^{n_c}$ are labeled.}
\label{fig:multiple_components}
\end{figure}

\subsection{Implementation Notes}
In this section, we describe the implementation details we used to produce the numerical results in ~Section~\ref{sec:results}. 
\subsubsection{Poles of Rational Curves--Implementation}
\label{sec:conversion}
In the test cases that follow, we consider $\mathbf{c}_i$ to be given in the Bernstein-\bezier\ basis. In order to find the roots of the control weight polynomial, \begin{equation*}
w_i(s) = \displaystyle\sum_{j=0}^{m_i} w_{i,j} B^{m_i}_j(s),
\end{equation*} 
we use Matlab's \texttt{roots} command~\cite{MATLAB}, which has an optimized algorithm for finding the eigenvalues of the companion matrix of $w_i(s)$. 
However, to use Matlab's \texttt{roots} command, we must convert from the Bernstein-\bezier\ basis to the monomial basis, which is an ill-conditioned basis conversion for high-degree polynomials, as discussed in~\cite{farouki1988algorithms}. The condition number of the conversion matrix is generally on the order of $10^{\frac{m}{2}}$, where $m$ is the degree of the polynomial. Approximately $\frac{m}{2}$ digits of precision are therefore lost in this conversion. The conversion matrix from the monomial basis to the Bernstein-\bezier\ basis is defined as
\begin{equation*}
MB_{i,j} = 
\begin{cases}\frac{\binom{n-k+1}{j-k}}{\binom{n}{j-1}} &\text{ for } j\geq k\\
0 &\text{ otherwise, }
\end{cases}
\end{equation*}
and the conversion from the Bernstein-\bezier\ basis to the monomial basis is defined as 
\begin{equation*}
BM_{i,j} = \begin{cases} 
(-1)^{j-k} \binom{n}{j-1} \binom{j-1}{k-1} & \text{ for } j\geq k\\
0 &\text{ otherwise, }
\end{cases}
\end{equation*}
where $BM_{i,j}$ denotes the entry in the $i$-th row and $j$-th column of the matrix. An optimal implementation would instead find roots in the original input basis to avoid basis conversion stability issues. 

\subsubsection{Green's Theorem--Implementation}
\label{sec:green_implement}
The choice of $C$ in the antiderivative computation in Equation~$\eqref{eq:antiderivative}$ affects the locations of the resulting quadrature points. If the $C$ is taken to be very far from the domain of interest, it can affect the floating point stability of the resulting computations, since the antiderivative integrals can become very large in magnitude. In order to keep the quadrature points close to the original domain, we choose $C$ to be the minimum y-component of any control point making up the curves in the boundary loop, $\Gamma$. That is, we choose
\begin{equation*}
C= \min_{i,j} y_{i,j}.
\end{equation*}  
The Gaussian quadrature points used to evaluate the integral in Equation~$\eqref{eq:antiderivative}$ will therefore always lie between the intermediate quadrature point $y$ and the smallest control point of the boundary loop $\Gamma$. Since the points $y$ always lie on the component curves, the quadrature points from the full algorithm will always lie within the bounding box of the control points. 

\subsubsection{Rational Intermediate Quadrature--Implementation}
\label{sec:rational_implementation}
To compute the intermediate quadrature rule for the modification to achieve polynomial exactness in Section~\ref{sec:rat_intermediate}, we use a method developed in \cite{deckers2017algorithm}. The method generates $m+1$ quadrature points when $m$ (possibly non-unique) poles are specified in linear time with respect to the number of poles $m$. The algorithm also ensures a degree of polynomial exactness beyond the default rational exactness. If $m$ poles are specified but $m+l$ quadrature points are requested, the algorithm will be exact for rational functions with the $m$ poles as well as for (one-dimensional) polynomial functions up to degree $l$. Convergence of the \spectralPE\ for arbitrary functions is achieved by increasing this $l$ parameter past that required by the degree of polynomial exactness $k$, for example in Figure~\ref{fig:nonpolynomial_results}.

The algorithm can be used to compute exact quadrature rules on any real interval for any class of rational functions with arbitrary poles outside of the domain of integration (in this case, the unit interval). It is well known that if all of the coefficients of a polynomial in the Bernstein-\bezier\ basis are of the same sign, then the polynomial has no roots within the unit interval \cite{farin2002curves}. Therefore, as long as all of the control weights $w_{i,j}$ in all of the component curves $\mathbf{c}_i$ are positive (a common situation in most applications), the quadrature scheme given in \cite{deckers2017algorithm} is theoretically valid.

\subsection{Remark: Location of Quadrature Points}
As can be seen in Figure~\ref{fig:quadrature_points_circle}, the quadrature points do not, in general, lie within the integration region $\Omega$. However, the quadrature points in our scheme are always contained within the bounding box of the component curves, $\tilde{\Omega}$. 
It is sometimes desirable that the locations of the quadrature points be contained within the boundaries of the quadrature region $\Omega$, particularly when dealing with potentially non-analytic integrands. If the original domain is convex, the locations of the quadrature points can be guaranteed to lie within the domain of integration using an appropriate choice of $C$ and a rotation of the coordinate system (see \cite{santin2011algebraic} for implementation details). If the domain is non-convex, a moment-fitting approach could be combined with our algorithm to compute a quadrature rule with only interior quadrature points \cite{thiagarajan2018shape, taber2018moment}. We briefly introduce moment-fitting and describe how it could be used with our method in \ref{app:momentFitting}.

\subsection{Remark: Multiple Components}
In some cases, the region of interest may not be genus-zero (i.e., it may have holes or disconnected components), see, for example, the region in Figure~\ref{fig:multiple_components}. In this case, there are multiple boundary loops $\partial \Omega = \{\Gamma_\alpha\}_{\alpha=1}^M$, each of which bounds a genus-zero domain $\Omega_\alpha$. The original domain $\Omega$ is then a finite Boolean (union, intersection, and subtraction) combination of the component genus-zero domains $\Omega_\alpha$. As long as the orientations of each of the boundary curves $\Gamma_\alpha$ are given, the algorithm as stated transfers seamlessly, since it treats each of the component curves $\mathbf{c}_i$ independently. In this case, the integral over $\Omega$ can  be thought of as the signed sum of the integrals over each $\Omega_\alpha$: $\iint_\Omega f(x,y) dxdy = \sum_{\alpha=1}^N s(\alpha) \iint_{\Omega_\alpha} f(x,y) dxdy$, where the indicator function $s(\alpha)$ can be sorted out according to the orientations of the boundary curves:
\begin{equation*}
s(\alpha) = \begin{cases} 1 & \text{ if $\Gamma_\alpha$ is counter-clockwise oriented}\\
-1 & \text{ if $\Gamma_\alpha$ is clockwise oriented}.
\end{cases}
\end{equation*} Our aforementioned algorithm is then applied to each integral individually.

\section{Numerical Test Cases}
\label{sec:results}

In this section, we test the \spectral\ and \spectralPE\ algorithms derived in Section~\ref{sec:algorithm} to compute integrals of various functions over six different domains. We first test the methods for monomial integrands on a circular region defined by four rational quadratic Bernstein-\bezier\ curves.  Then, we investigate the methods on five other relatively more complicated regions: 
\begin{inparaenum}[(1)]
	\item a plate with a hole, commonly used as an isogeometric test case \cite{hughes2005isogeometric}, 
	\item an L-bracket with three holes, also used as an isogeomtric test case \cite{benzaken2017rapid}, 
	\item a guitar-shaped region,
	\item a treble clef-shaped region, and 
	\item a region representing a cross section of two interlocked screws \cite{hinz2018spline}.
\end{inparaenum} In our experiments, our quadrature schemes display the same properties on simpler shapes, such as the region in Figure~\ref{fig:immersogeometric_example}, as on the six test shapes in this section.

We show numerically that our methods produce spectrally convergent quadrature schemes. In addition we show that the \spectralPE\ method, with exact rational intermediate quadrature, is exact for polynomials up to degree $k$ when
\begin{equation}
  n_q = \left\lceil{\frac{k+1}{2}}\right\rceil\sum_{i=1}^{n_c} \left(m_i(k + 3)+1\right) \label{eq:n_q}
\end{equation}
quadrature points are used for regions enclosed by $n_c$ curves of degree $m_i$ and maximum polynomial integrand degree $k$. 
 
 In the numerical tests, including for the results displayed in Figures~\ref{fig:quadrature_comparison_circle},~\ref{fig:quadrature_weird_shapes},~\ref{fig:nonpolynomial_results},~and~\ref{fig:nonpolynomial_timing}, we refer to our \spectral\ scheme as \textsc{GT-}\spectral\ and our \spectralPE\ scheme as \textsc{GT-}\spectralPE, since they rely on Green's theorem (GT). We define integration error as the absolute value of the difference between the correct integral and the calculated integral. We calculate the correct value of each integral as the value obtained when using $P=Q_i=55$ quadrature points in each component quadrature rule with our \spectral\ algorithm. All plots in this section are Log-Log plots to emphasize the orders of magnitude differences in accuracy and efficiency among the various methods. 
 
For polynomial integrands with our \spectralPE\ method, we only calculate one quadrature rule which uses the number of quadrature points required to integrate the given integrand exactly up to machine precision, given by the formula in Equation~$\eqref{eq:n_q}$. For all integrands with our \spectral\ method and for nonpolynomial integrands with our \spectralPE\ method, we increase the number of quadrature points in each component quadrature rule according to the strategies given in Sections~\ref{sec:intermediate},~\ref{sec:antiderivative},~and~\ref{sec:rational_implementation}. In particular, we let the number of quadrature points in the antiderivative and intermediate quadrature rules be equal and gradually increase both in order to obtain convergence.

\subsection{Comparison Methods} 
\label{sec:comparison}
To show the relative benefit of our methods, we compare with three domain decomposition-based (\textsc{DD}) quadrature schemes and two Green's theorem-based (\textsc{GT}) quadrature schemes from the literature, which represent the state of the art for integrating over rational regions.
\begin{compactdesc}
  \item[\textsc{DD-Quadtree}.] Our implementation of quadtree integration refines a background Cartesian grid uniformly, then performs three levels of adaptive refinement anywhere where the boundary intersects the background grid and uses 3$^{rd}$-order (4 point) Gaussian quadrature on each element~\cite{kudela2015efficient}.
  \item[\textsc{DD-Linear mesh}.] We use an implementation of linear meshing from \emph{mesh2d} which refines the discretization by remeshing the geometry with a specified maximum element size and uses 3$^{rd}$-order (4 point) Gaussian quadrature on each element~\cite{engwirda2014locally}.
  \item[\textsc{DD-Rational mesh}.] We use an implementation of exact rational meshing from TRIGA with increasing order Gaussian quadrature on each element~\cite{engvall2016isogeometric}.
  \item[\textsc{GT-Linear}.] To compare with a Green's theorem-based method with linear geometric approximations, we use \emph{SplineGauss}~\cite{sommariva2009gauss}, which approximates the boundary with linear segments and uses a high-order Gaussian quadrature rule to evaluate the antiderivatives \cite{sommariva2009gauss}. In their terminology, we use $5^{th}$ degree of precision for the polynomial integrands and $15^{th}$ degree of precision for the nonpolynomial integrands.
  \item[\textsc{GT-Cubic Spline}.] To compare with a Green's theorem-based method with high-order geometric approximations, we also use \emph{SplineGauss}, but with a cubic spline approximation of the boundary instead of linear~\cite{sommariva2009gauss}. For this method, we use 10$^{th}$ degree of precision for the polynomial integrands and 22$^{nd}$ degree of precision for the nonpolynomial integrands.
\end{compactdesc}
   We note that we have not optimized any of the implementations of these algorithm (nor have we optimized implementations of our own algorithms), but that we expect the broad trends, particularly with respect to differences in pre-processing time and convergence rates, to hold true in any implementation. In our experiments, using higher-order Gaussian quadrature rules for the \textsc{DD-Quadtree}, \textsc{DD-Linear mesh}, \textsc{GT-Cubic Spline}, and \textsc{GT-Linear} strategies did not improve the error much, but yielded a much higher number of quadrature points. This implies that the geometric approximation itself is the primary source of error in these methods.
 
\subsection{Exact Monomial Quadrature over a Circle}
\label{sec:results_circle}

First, we test the algorithm in-depth on a simply connected, convex test case to show its fundamental properties. We integrate the monomial functions $1, x, y , x^2, \ldots$ up to degree $5$ over a circular region $\Omega$ centered at $(x_0,y_0)$ with radius 1. See \ref{app:detailed_example} for an in-depth, detailed work-through for this test case.

\subsubsection{Domain Definition} 
\label{sec:circle_def}
The unit circle is defined by $n_c=4$ quadratic rational Bernstein-\bezier\ curves $\mathbf{c}_0,\mathbf{c}_1\mathbf{c}_2,\mathbf{c}_3$ such that $c_i = c_i' - (x_0,y_0)$ where
\begin{align*}
\mathbf{c}'_0 = \begin{cases}\left( \frac{B_0(s) + \frac{\sqrt{2}}{2} B_1(s)}{B_0(s) + \frac{\sqrt{2}}{2} B_1(s) + B_2(s)} \right)\\
\left(\frac{\frac{\sqrt{2}}{2} B_1(s) + B_2(s)}{B_0(s) + \frac{\sqrt{2}}{2} B_1(s) + B_2(s)}\right)\end{cases} \mathbf{c}'_1 = \begin{cases} \left(\frac{-\frac{\sqrt{2}}{2} B_1(s) - B_2(s)}{B_0(s) + \frac{\sqrt{2}}{2} B_1(s) + B_2(s)}\right)\\
\left( \frac{B_0(s) + \frac{\sqrt{2}}{2} B_1(s)}{B_0(s) + \frac{\sqrt{2}}{2} B_1(s) + B_2(s)}\right)\end{cases}\\
 \mathbf{c}'_2 = \begin{cases} \left(\frac{-B_0(s) - \frac{\sqrt{2}}{2} B_1(s)}{B_0(s) + \frac{\sqrt{2}}{2} B_1(s) + B_2(s)}\right)\\
\left(\frac{-\frac{\sqrt{2}}{2} B_1(s) - B_2(s)}{B_0(s) + \frac{\sqrt{2}}{2} B_1(s) + B_2(s)}\right)\end{cases} \mathbf{c}'_3 = \begin{cases}\left(\frac{\frac{\sqrt{2}}{2} B_1(s) + B_2(s)}{B_0(s) + \frac{\sqrt{2}}{2} B_1(s) + B_2(s)}\right)\\
\left( \frac{-B_0(s) - \frac{\sqrt{2}}{2} B_1(s)}{B_0(s) + \frac{\sqrt{2}}{2} B_1(s) + B_2(s)}\right),\end{cases}\\
\end{align*}
and $B_0(s) = (1-s)^2$, $B_1(s) = 2(1-s)s$, and $B_2(s) = s^2$. 

\subsubsection{Results}
\begin{figure}

\begin{center}
\includegraphics[width=.75\linewidth]{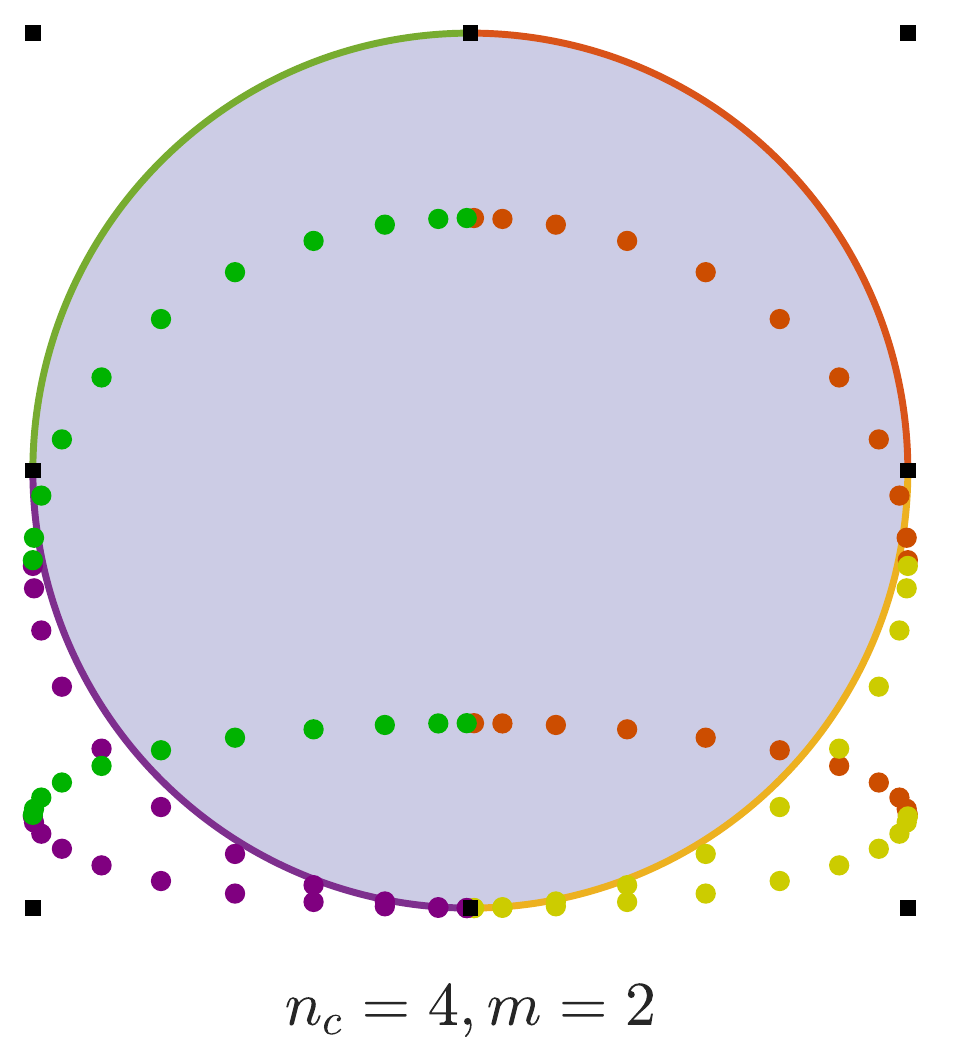}
\end{center}
\caption{\spectralPE\ quadrature points for a circle defined by 4 rational parametric Bernstein-\bezier\ curves integrating polynomials up to degree $k=3$ exactly. Quadrature points are contained within the bounding box of the curve control points (black squares). The quadrature points are color-coded according to their corresponding component curve in the evaluation of Green's theorem. } 
\label{fig:quadrature_points_circle}
\centering
\vspace{.4cm}
\includegraphics[width=.95\linewidth]{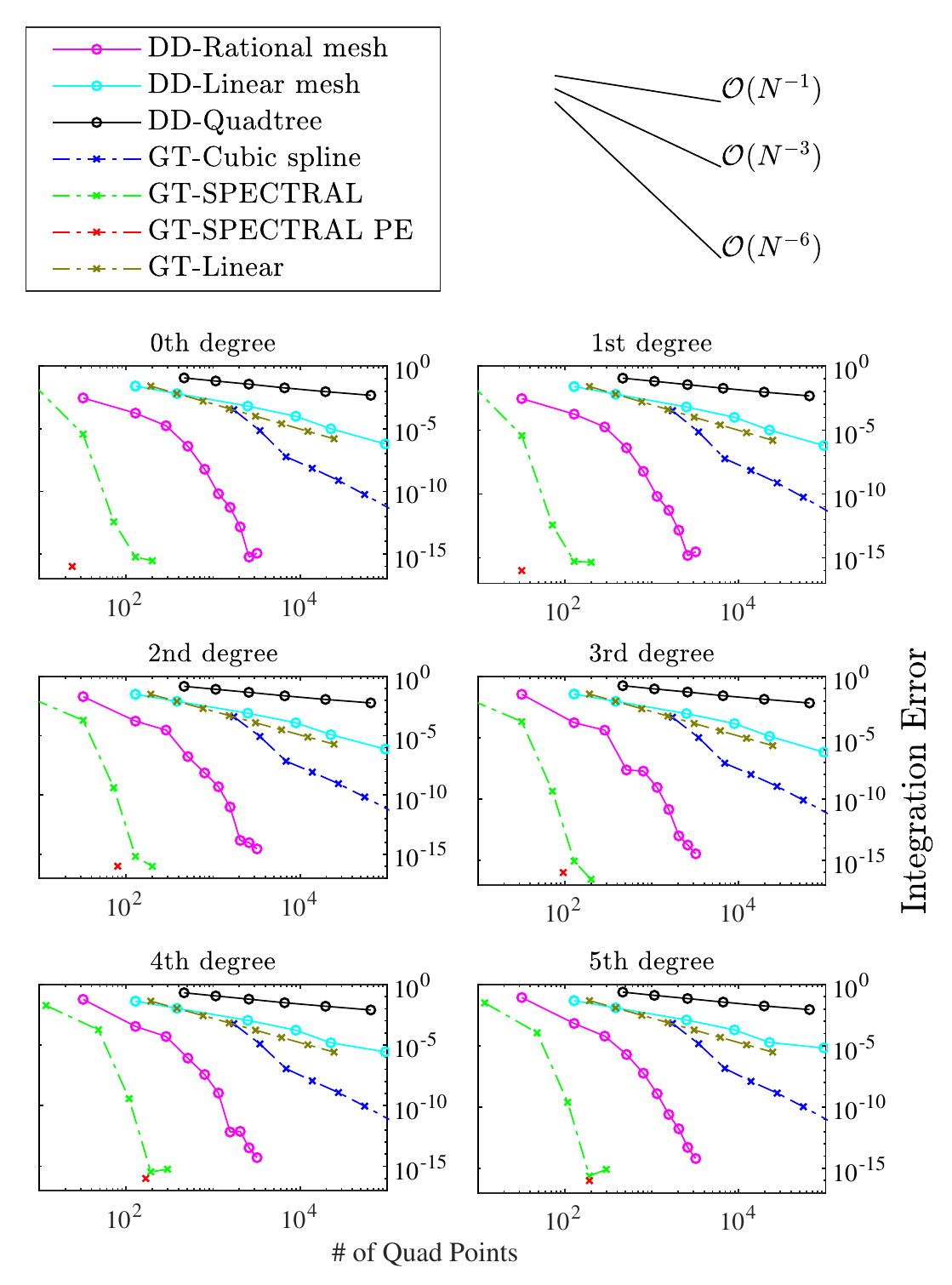}
\vspace{-.3cm}
\caption{Integrals of monomials over the circle can be computed exactly using the \spectralPE\ algorithm. Moreover, the \spectral\ algorithm achieves faster convergence than any algebraic order and is much more efficient than other common methods. Note that the error values are averages over the monomials of a particular degree. In this case, all curves are of degree $m_i=2$.}
\label{fig:quadrature_comparison_circle}
\end{figure}

The quadrature points produced using the \spectralPE\ algorithm with a maximal degree of polynomial exactness of $k=3$ are shown in Figure~\ref{fig:quadrature_points_circle}. The \spectral\ algorithm produces similarly distributed quadrature points.


As can be seen in Figure~\ref{fig:quadrature_comparison_circle}, the \spectralPE\ quadrature strategy yields exact results with the expected number of quadrature points 
$n_q = 28$ quadrature points in the case of $p(x,y) =1$ and $n_q = 36$ quadrature points in the case of $p(x,y) =x,y$, as per Eq.~$\eqref{eq:n_q}$. Moreover, the \spectral\ algorithm produces quadrature rules which converge faster than any algebraic order with respect to number of quadrature points and are orders of magnitude more efficient than the comparison methods.

\subsection{Exact Polynomial Quadrature over Rational Shapes}
\label{sec:results_rational_shapes}
We also test our algorithms on a variety of more complicated shapes, including
\begin{inparaenum}[(1)]
	\item a plate with a hole, 
	\item an L-bracket with three holes, 
	\item a guitar-shaped region,
	\item a treble clef-shaped region, and
	\item a region representing a cross section of two interlocked screws~\cite{hinz2018spline}.
\end{inparaenum}
The regions are displayed in Figure~\ref{fig:quadrature_weird_shapes}.

\subsubsection{Domain definition}
The first two domains (plate with hole and L-bracket with 3 holes) were formed by taking finite Boolean (i.e.\ unions, intersections, and differences) combinations of quadratic rational \bezier\ regions with known control weights and control points (namely, rectangles and ellipses). Booleans of arbitrary rational Bernstein-\bezier\ shapes can be computed using methods briefly outlined in~\ref{app:booleans}. This same procedure could be used to find boundary parameterizations for intersections of elements in, for example, immersogeometric analysis. The plate with a hole is composed of $M=2$ boundary loops and the boundary contains a total of $n_c=8$ rational \bezier\ curves with maximum degree $m=2$. The L-bracket with $3$ holes is composed of $M=4$ boundary loops and the boundary contains a total of $n_c=20$ rational \bezier\ curves with maximum degree $m=2$. 


The guitar, treble clef, and screws regions were designed by editing the control weights of NURBS curves, then projecting the NURBS curves onto rational \bezier\ curves. The algorithm for converting NURBS to rational Bernstein-\bezier\ curves is briefly discussed in \ref{app:NURBS}. The guitar region is composed of $M=2$ boundary loops with a total of $n_c= 60$ component  rational \bezier\ curves with maximum degree $m=3$. The treble clef region is composed of $M=4$ boundary loops containing $n_c=99$ rational \bezier\ curves with maximum degree $m=2$. The two interlocked screws region is composed of $M=3$ boundary loops with a total of $n_c=128$ rational \bezier\ curves with maximum degree $m=3$. 


\subsubsection{Results}
\label{sec:results_rational_shapes_results}

\begin{figure*}
\vspace{-1.5cm}
\includegraphics[width=1\linewidth]{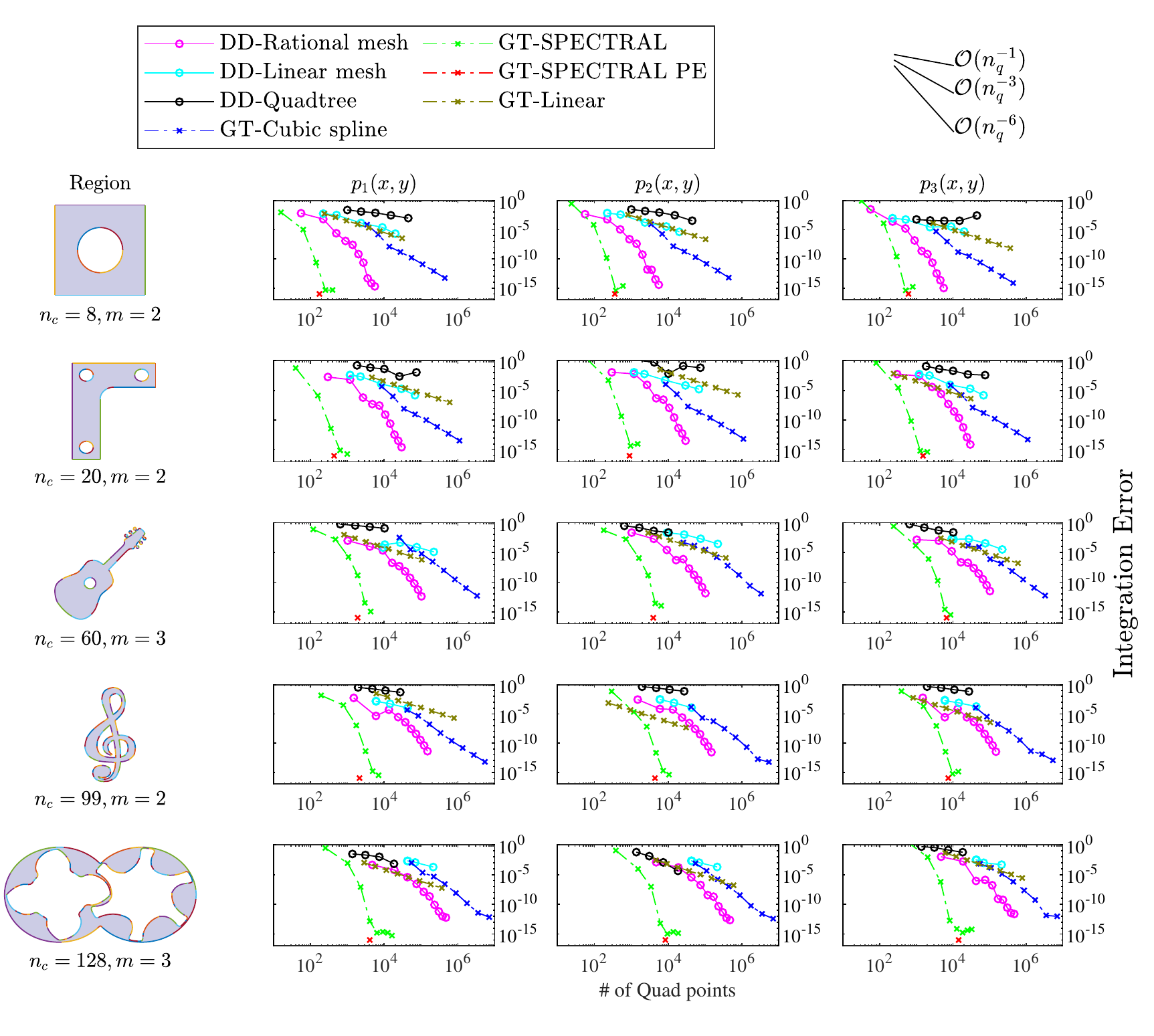}
\caption{Integrals of polynomials over arbitrary rational parametric regions can be computed exactly using the presented \spectralPE\ quadrature scheme. Moreover, the \spectral\ scheme converges faster than any algebraic order and is much more efficient than other methods in the literature. The left column shows the region being integrated in light grey and the different colors on the boundaries represent the component curves $\{\mathbf{c}_i\}_{i=1}^{n_c}$. See Section~\ref{sec:comparison} for a more detailed description of the comparison methods. The integrands are $p_1(x,y)=2x^2 + xy -y +2$, $p_2(x,y)=2x^2y^2+\frac{3}{10}x^2y-y^4+3x+2$, and $p_3(x,y)=x^5-5y^3x^3+2yx^2+\frac{1}{5}x^2+3$.}
\label{fig:quadrature_weird_shapes}
\end{figure*}

We applied the \spectral\ and \spectralPE\ algorithms 
to the shapes defined above for various polynomials of degrees $2$, $4$, and $6$. As can be seen in Figure~\ref{fig:quadrature_weird_shapes}, both of the presented quadrature schemes outperform the five comparison methods in terms of number of quadrature points for a given degree of accuracy by orders of magnitude. In addition, the \spectralPE\ scheme is able to achieve exactness with the expected number of quadrature points. We observed the same behavior for a variety of polynomial integrands, including higher-degree polynomials. We note here that the comparison methods achieve the expected convergence rates of 1$^{st}$ order convergence for \textsc{DD-Linear mesh}, \textsc{DD-Quadtree}, and \textsc{GT-Linear}, 3$^{rd}$ order convergence for \textsc{GT-Cubic spline}, and spectral convergence for \textsc{GT-Rational mesh}.

\subsection{Quadrature of Non-polynomial Integrands over Rational Shapes}
\label{sec:results_non_polynomial_integrands}

We also compare our algorithms on a representative set of non-polynomial integrands, including a rational integrand, an exponential integrand, and a square root integrand. These integrands require adaptive integration, even for the \spectralPE\ algorithm. The \textit{a priori} formula for the number of quadrature points will not guarantee exactness, since the antiderivative computation using Gaussian quadrature is no longer exact in this case. We therefore use varying amounts of quadrature points in the \spectralPE\ algorithm to observe convergence of the method. To do this, we find a quadrature rule exact for polynomials up to degree $k=2$, then increase the number of quadrature points in the intermediate rational quadrature algorithm as explained in Section~\ref{sec:rational_implementation} (by increasing $l$) and in the antiderivative Gaussian quadrature order as explained in Section~\ref{sec:antiderivative} (by increasing $P$), while keeping the number of quadrature points in each component rule equal.

\begin{figure*}
\vspace{-1.5cm}
\includegraphics[width=1\linewidth]{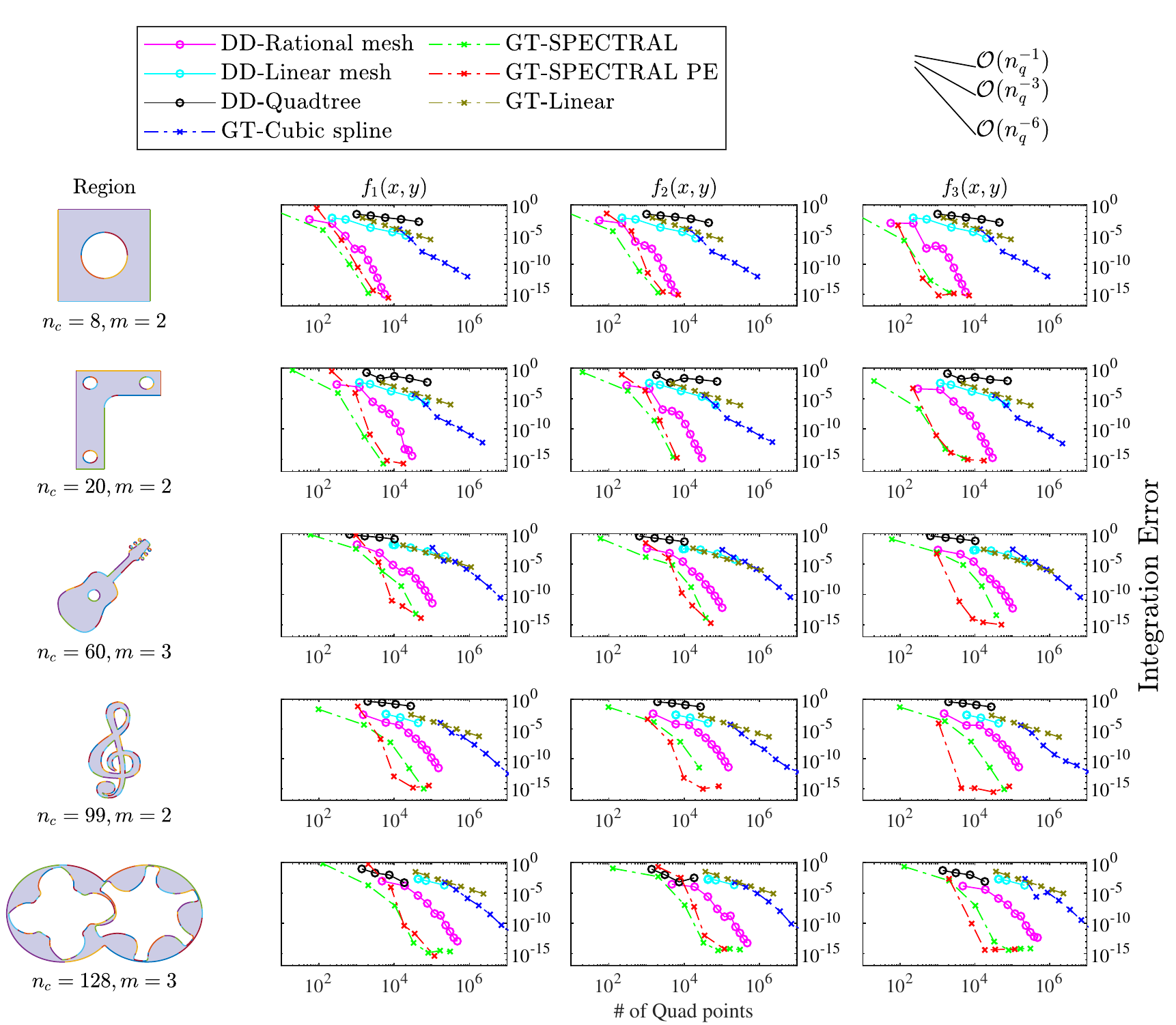}
\caption{Both the \spectral\ and the \spectralPE\, methods given in Section~\ref{sec:algorithm} converge faster than any algebraic order on a broad class of nonpolynomial integrands. The left column shows the region being integrated in light grey and the different colors on the boundaries represent the component curves $\{\mathbf{c}_i\}_{i=1}^{n_c}$. See Section~\ref{sec:comparison} for a more detailed description of the comparison methods. The integrands are $f_1(x,y)=\frac{y^3-x^3y^2-xy-3}{x^2y^2+30}$, $f_2(x,y)=e^{-x^2+2y}$, and $f_3(x,y)=\sqrt{(x+10)^2+(x+10)(y+10) +x}$.}  
\label{fig:nonpolynomial_results}
\end{figure*}

As can be seen in Figure~\ref{fig:nonpolynomial_results}, both the \spectral\ and \spectralPE\ methods are orders of magnitude more efficient than other methods in the literature on a per-quadrature point basis and converge to machine precision with far fewer quadrature points. We observed similar behavior for a variety of non-polynomial integrands. In this case, it is likely more beneficial to use the \spectral\ algorithm, since it has the same convergence behavior as the \spectralPE\ algorithm, but is likely relatively simpler to implement.

\subsection{Timing Results}
\begin{figure}
\includegraphics[width=1\linewidth]{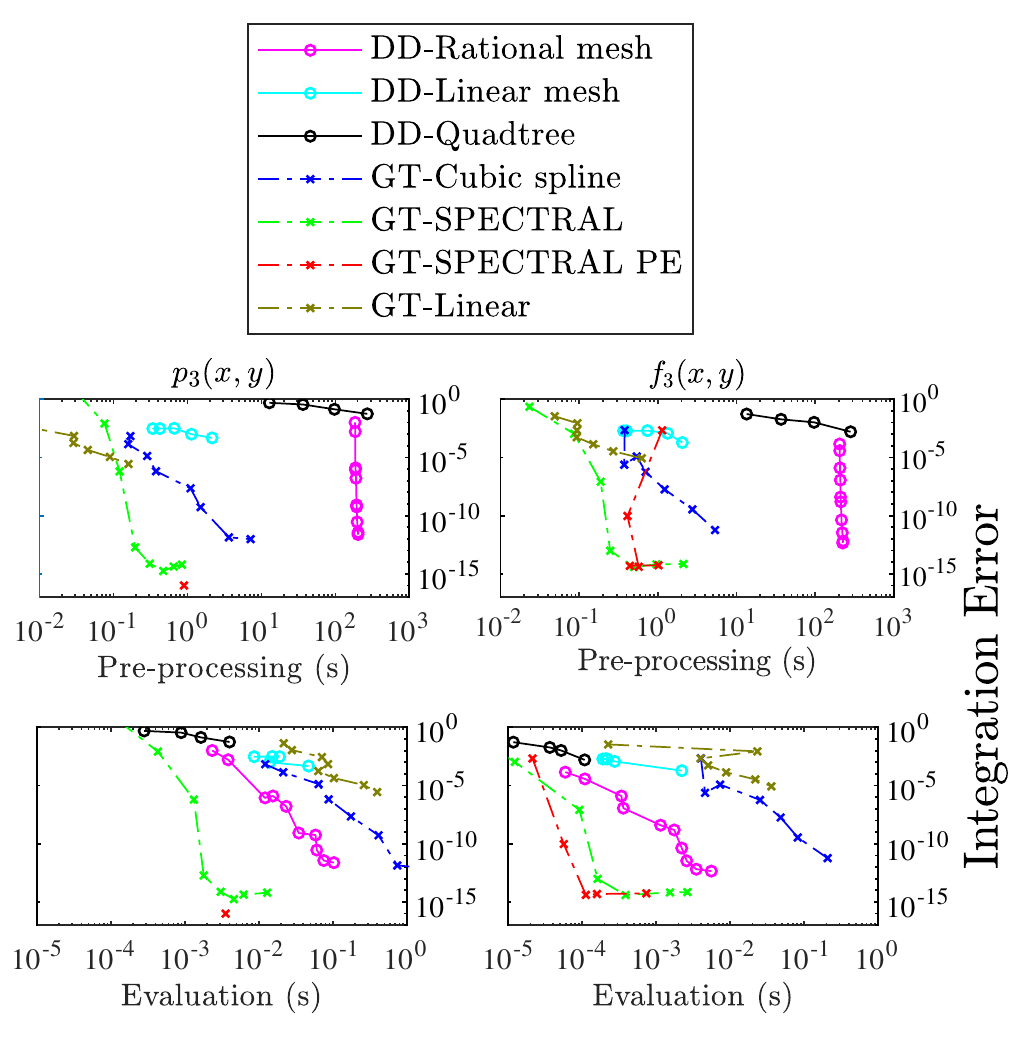}
\caption{Our proposed methods are much more efficient for a given error in terms of computational time. Pre-processing (top) and evaluation (bottom) timings to integrate functions $p_3$ and $f_3$ for the interlocked screws mesh. Timings correspond to the lower-right most plots of Figure~\ref{fig:quadrature_weird_shapes} (left) and Figure~\ref{fig:nonpolynomial_results} (right). }
\label{fig:nonpolynomial_timing}
\end{figure}

In addition to performing more efficiently with respect to the number of quadrature points, we found that both the \spectral\ and \spectralPE\ quadrature schemes are also orders of magnitude more efficient with respect to computational time required for nearly all cases shown in the previous sections. As representative examples, we show timing results in Figure~\ref{fig:nonpolynomial_timing} to calculate the various quadrature points and weights (pre-processing) and apply them to a single function (evaluation) for the lower-rightmost plot in each of Figures~\ref{fig:quadrature_weird_shapes} and \ref{fig:nonpolynomial_results}. Timings were performed using MATLAB's \lstinline{timeit} function. 
As can be seen, the present algorithms require less evaluation time than any of the comparison methods for a given level of accuracy, except for \textsc{GT-Linear} for very low accuracies. In addition, the present algorithms require far less pre-processing time than the only other spectral method, rational meshing. 

Timing tests were performed using MATLAB version R2019a on Windows 10 with 32GB of RAM. It is important to note that the timings given here use non-optimized MATLAB code. Results would likely change with optimized compiled implementations, although we don't expect changes in the trends or relative performance of the various algorithms, particularly with respect to convergence behavior and pre-processing time.

\section{Discussion}
\label{sec:discussion}

In this work, we have presented two spectral meshless algorithms for computing quadrature rules for arbitrary planar regions bounded by rational curves. Both the \spectral\ and \spectralPE\ algorithms converge faster than any algebraic order with respect to total number of quadrature points for arbitrary smooth integrals. Moreover, the \spectralPE\ algorithm is exact (up to machine precision) for polynomial integrands up to a pre-specified degree, $k$. 
In contrast to other schemes in the literature, which either require expensive meshing to achieve spectral convergence or have algebraic orders of convergence, 
both of our algorithms achieve spectral convergence without meshing. Moreover, our \spectralPE\ method prescribes the number of quadrature points necessary to exactly integrate polynomial functions over these regions, whereas other methods must spend extra effort adaptively converging to the correct solution.

We presented numerical test cases which show that the \spectralPE\ algorithm is exact for computing integrals of polynomials over rational geometries, reducing the number of quadrature points necessary to achieve a given precision for polynomial integrands even compared with our \spectral\ rule. We also presented numerical test cases which show that both the \spectral\ and \spectralPE\ algorithms presented here are much more efficient than other methods when integrating arbitrary functions over rational geometry, both on a per-quadrature point basis and a computational time basis. 

Our intention is for this paper to be the first in a series which describes approaches for high-order or spectrally convergent numerical quadrature over general parametric regions in $\mathbb{R}^2$ and $\mathbb{R}^3$. The broad idea is to iteratively use the parametric version of generalized Stokes' theorem combined with numerical antiderivative computation to compute efficient quadrature rules. This paper outlines the approach applied to integration over arbitrary two-dimensional rational (and, as a special case, polynomial) parametric regions and achieves spectral convergence. Part two will address the general three-dimensional parametric case using an extension of the algorithm presented here and we expect it to attain high algebraic orders of convergence.

Another potential avenue for future research could include expanding the \spectralPE\ scheme for polynomial integrands to 3D rational geometries, although the poles of rational parametric surfaces can in general be algebraic curves, for which exact rational quadrature has not been studied extensively.
We would also like to apply these methods to immersogeometric analysis, field transfer between high-order meshes, 
and initialization of material volume fractions in multimaterial simulations.

\section*{Acknowledgments}
This work was performed under the auspices of the U.S. Department of Energy by Lawrence Livermore National Laboratory under Contract DE-AC52-07NA27344. 

We would like to thank the reviewers for their thoughtful comments and suggestions in the preparation of this manuscript. We also thank the authors of the open source Matlab software that was used for our \spectralPE\ intermediate quadrature rule and for the comparison methods, including \textit{rfejer}~\cite{deckers2017algorithm}, \textit{mesh2d}~\cite{engwirda2014locally}, TRIGA~\cite{engvall2016isogeometric}, and \textit{SplineGauss}~\cite{sommariva2009gauss}.

\appendix
\section{Booleans of Rational Shapes}
\label{app:booleans}
Intersections, unions, and differences of rational shapes can be computed using the algorithm given in Section~\ref{sec:algorithm} by using \bezier\ intersection and reparametrization. 
An algorithm for computing parametrizations of Boolean combinations (i.e. unions, intersections, and differences) of regions defined by polynomial triangular Bernstein-\bezier\ elements is described in~\cite{hermes2018high}. We have adapted this algorithm to the rational case to help design some of the shapes for our empirical results in Section~\ref{sec:results}. Our approach combined with efficient intersection could be used to efficiently perform the integration necessary in immersogeometric analysis, field transfer between high-order meshes, and VOF field initialization, in which integrals over intersections between rational geometries must be computed.
\section{NURBS}
\label{app:NURBS}

The algorithm described in Section~\ref{sec:algorithm} can easily be used on regions bounded by non-uniform rational basis splines (NURBS). Because NURBS can be decomposed exactly into rational \bezier\ curves, the algorithm only requires one extra step: \bezier\ extraction~\cite{borden2011isogeometric}, which converts a degree-$m$ NURBS curve into component degree-$m$ rational \bezier\ curves.  

\section{Moment Fitting}
\label{app:momentFitting}
The \spectral\ and \spectralPE\ algorithms described in Section~\ref{sec:algorithm} can also be used as an input to moment fitting, through which interiority of the quadrature points can be achieved. Moment fitting is a technique that was reintroduced to the engineering community in \cite{mousavi2010generalized}. In moment fitting, a linear system of equations of the form
\begin{equation*}
\begin{bmatrix} 1 & 1 & 1 &\ldots & 1\\
 x_1 & x_2 & x_3 &\ldots & x_n \\
 y_1 & y_2 & y_3 &\ldots & y_n \\
 x_1^2 & x_2^2 & x_3^2 &\ldots & x_n^2\\
 x_1y_1 & x_2y_2 & x_3y_3 & \ldots &  x_ny_n\\ 
 y_1^2 & y_2^2 & y_3^2& \ldots & y_n^2\\
\vdots & \vdots & \vdots & \ddots & \vdots \\
\end{bmatrix}\begin{bmatrix} w_1 \\ w_2\\ w_3 \\\vdots \\ w_n \end{bmatrix} = \begin{bmatrix}
\int_D 1 dA\\
\int_D x dA\\
\int_D y dA\\
\int_D x^2 dA\\
\int_D xy dA\\
\int_D y^2 dA,\\
\vdots\\
\end{bmatrix}
\end{equation*}
is used to find corresponding quadrature weights $\{w_i\}_{i=1}^n$ for a given pre-specified set of quadrature points $\{ (x_i,y_i)\}_{i=1}^n$ for integrating over a region $D$. The number of quadrature points $n$ is typically taken to be much larger than the number of monomials fitted. The quadrature points $\{ (x_i,y_i)\}_{i=1}^n$ can be specified to be within the domain of integration $D$ \cite{thiagarajan2018shape}. However, the quadrature points should be chosen with care, as their locations impact the accuracy and stability of the resulting quadrature scheme for general integrands. Note that non-negative quadrature weights can be ensured through the solution of a related constrained optimization problem \cite{keshavarzzadeh2018numerical}.

Our method can obtain exact results for the right-hand side vector for rational domains $D$, eliminating any error in the moment-fitting procedure stemming from monomial integration and making it a prime candidate for evaluating the right-hand side. Moreover, the lack of a guarantee of interiority of quadrature points in our method is one of its potential disadvantages, so using it with moment-fitting can be helpful in contexts when interiority is needed. We leave exploration of such a method to a future work.
\section{Detailed, Step-by-step Application of the \spectralPE\ Algorithm over a Circle}
\label{app:detailed_example}
In this section, we describe the steps of the \spectralPE\ algorithm in detail for the integral $\int_\Omega xy^2 dxdy$, where $\Omega$ is the circular region defined in Section~\ref{sec:circle_def} with center $(x_0,y_0) = (0,0)$ and radius $1$. The steps in this section are visualized in Figure~\ref{fig:overview} and the final quadrature rule produced is shown in Figure~\ref{fig:quadrature_points_circle}. Note that the maximum polynomial degree of the integrand is $k=3$. As described in Section~\ref{sec:algorithm}, the \spectralPE\ algorithm for computing polynomially-exact quadrature rules proceeds for each component curve $\mathbf{c}_i$ in two steps: \begin{enumerate}[(1)]
  \item Compute a set of intermediate quadrature points and weights for the Green's Theorem line integral corresponding to $\mathbf{c}_i$ and
  \item Compute quadrature rules to exactly evaluate the antiderivative function numerically at each intermediate quadrature point.
\end{enumerate} 

Let us first consider the curve $\mathbf{c}_0$, noting that the following procedure must be applied to each of the $\{\mathbf{c}_i\}_{i=0}^3$. We have
\begin{equation*}
\mathbf{c}_0 = \begin{cases} x_0(s) = \left( \frac{(1-s)^2 + \sqrt{2} (1-s)s}{(1-s)^2 +  \sqrt{2} (1-s)s + s^2} \right) = \left( \frac{1+(\sqrt{2}-2)s + (1-\sqrt{2})s^2}{ 1+(\sqrt{2}-2)s + (2-\sqrt{2})s^2} \right)\\
y_0(s) = \left(\frac{\sqrt{2} (1-s)s + s^2}{(1-s)^2 +  \sqrt{2} (1-s)s + s^2}\right) = \left( \frac{\sqrt{2}s + (1-\sqrt{2})s^2}{1+(\sqrt{2}-2)s + (2-\sqrt{2})s^2} \right).\end{cases}
\end{equation*} For the \spectralPE\ algorithm, the intermediate quadrature rule must be exact for the intermediate rational function appearing in the Green's theorem line integral in Equation~$\eqref{eq:green_par}$. In the actual algorithm, the antiderivative and the partial derivative terms are computed numerically. However, for demonstration purposes, we show the Green's theorem line integral for this curve and integrand with analytic antiderivative and partial derivative terms:
\begin{align*}
&\int_{\mathbf{c}_i} A_f(x(s),y(s))\frac{d x_i}{ds} ds = \int_0^1 \frac{1}{3}x_i(s) y_i(s) ^3 \frac{d x_i(s)}{ds} ds\\
= \int_0^1 \frac{1}{3} &\left( \frac{p(s)}{\left( 1+(\sqrt{2}-2)s + (2-\sqrt{2})s^2\right)^6} \right)ds,
\end{align*}
where $p(s)$ is a polynomial of degree $10$ which is too long to write out here, but is not important. Importantly, the integrand function has the same poles as the original Bernstein-\bezier\ curve, but with multiplicities multiplied by $k+3=6$. The next step in the algorithm is therefore to find these poles. In general, we use a numerical routine described in Sections~\ref{subsec:poles}~and~\ref{sec:conversion}. In this explanatory example, we find them analytically for demonstration purposes as
\begin{equation*}
p_{0,0} = \frac{1}{2} +\frac{\sqrt{2}}{4- 2\sqrt{2}}i,\hspace{.5cm}  p_{0,1} = \frac{1}{2} -\frac{\sqrt{2}}{4- 2\sqrt{2}}i.
\end{equation*} 
Finally, the rational functions which must be integrated exactly are those functions with poles $p_{0,0}$, $p_{0,1}$ each having multiplicity up to $6$. In general, we use the \textit{rfejer} routine to calculate the quadrature nodes and weights \cite{deckers2017algorithm}. In this case, there will be $13$ quadrature nodes and weights which exactly integrate this class of rational functions.

In the above derivation, we computed the partial derivative and antiderivative terms analytically. In general, these are not known analytically. Therefore, we generally numerically evaluate each of these as needed for each intermediate quadrature point. For example, the first intermediate quadrature point for $\mathbf{c}_0$ is $s_{0,0}= 0.004210269296207$ with weight $\gamma_{0,0} = 0.014622491933143$. In order to find the full quadrature rule, we must evaluate both the antiderivative and derivative terms at this point. 

For the antiderivative computation, we pick $C=-1$, because the lowest control point defining the entire domain has $y$-coordinate $-1$:
\begin{equation*}
A_f(x(s_{0,0}),y(s_{0,0}))=\int_{-1}^{y(s_{0,0})} f(x(s_{0,0}),t)dt.
\end{equation*}
This will cause the resulting final quadrature points for this intermediate quadrature point to be spread between $y=-1$ and $y=y(s_{0,0})$. To further clarify the effect of picking a particular $C$ value on the locations of the quadrature points, compare the quadrature points for the orange upper-right curve in Figure~\ref{fig:overview}, for which we took $C=-1$, to the quadrature points in Figure~\ref{fig:quadrature_points_circle}, in which we took $C=0$.	

  Because we know $f(x,y)$ is a polynomial of degree $3$, the number of Gaussian quadrature points necessary for exactness is $2$. The points and weights for the antiderivative quadrature rule for $s_{0,0}$ are
\begin{align*}
x_{0,0,0} &= \phantom{-}0.005961518538852,\hspace{.5cm}  x_{0,0,1} = \phantom{-}0.005961518538852,\\
y_{0,0,0} &= -0.577354024434498, \hspace{.5cm} y_{0,0,1} = \phantom{-}0.577336254424967,\\
\gamma_{0,0,0} &= \phantom{-}0.502980759269426, \hspace{.5cm} \gamma_{0,0,1} = \phantom{-}	0.502980759269426.
\end{align*}

The derivative at $s_{0,0}=0.008451616990154$ can also be evaluated numerically via deCasteljau's algorithm as $\frac{dx_0(s_{0,0})}{ds} = 1.417670137190716$. Therefore, the final quadrature weights for the two quadrature points, $(x_{0,0,0}, y_{0,0,0})$ and $(x_{0,0,1}, y_{0,0,1})$ are, respectively,  
\begin{align*}
w_{0,0,0} &= \gamma_{0,0}\gamma_{0,0,0}\frac{dx_0(s_{0,0})}{ds} = 0.00006216022389335293,\\
w_{0,0,1} &= \gamma_{0,0}\gamma_{0,0,1}\frac{dx_0(s_{0,0})}{ds} = 0.00006216022389335293.
\end{align*} 

This process would be repeated for the other intermediate quadrature points and other curves. 

\bibliography{refs}

\begin{thebibliography}{10}
\expandafter\ifx\csname url\endcsname\relax
  \def\url#1{\texttt{#1}}\fi
\expandafter\ifx\csname urlprefix\endcsname\relax\def\urlprefix{URL }\fi
\expandafter\ifx\csname href\endcsname\relax
  \def\href#1#2{#2} \def\path#1{#1}\fi

\bibitem{hughes2005isogeometric}
T.~J. Hughes, J.~A. Cottrell, Y.~Bazilevs, Isogeometric analysis: {CAD}, finite
  elements, {NURBS}, exact geometry and mesh refinement, {Computer Methods in
  Applied Mechanics and Engineering} 194~(39-41) (2005) 4135--4195.

\bibitem{bazilevs2010isogeometric}
Y.~Bazilevs, V.~M. Calo, J.~A. Cottrell, J.~A. Evans, T.~J.~R. Hughes,
  S.~Lipton, M.~A. Scott, T.~W. Sederberg, Isogeometric analysis using
  {T}-splines, {Computer Methods in Applied Mechanics and Engineering}
  199~(5-8) (2010) 229--263.

\bibitem{borden2011isogeometric}
M.~J. Borden, M.~A. Scott, J.~A. Evans, T.~J. Hughes, Isogeometric finite
  element data structures based on {B{\'e}zier} extraction of {NURBS},
  International Journal for Numerical Methods in Engineering 87~(1-5) (2011)
  15--47.

\bibitem{kamensky2015immersogeometric}
D.~Kamensky, M.-C. Hsu, D.~Schillinger, J.~A. Evans, A.~Aggarwal, Y.~Bazilevs,
  M.~S. Sacks, T.~J. Hughes, An immersogeometric variational framework for
  fluid--structure interaction: Application to bioprosthetic heart valves,
  {Computer Methods in Applied Mechanics and Engineering} 284 (2015)
  1005--1053.

\bibitem{schillinger2012isogeometric}
D.~Schillinger, L.~Dede, M.~A. Scott, J.~A. Evans, M.~J. Borden, E.~Rank, T.~J.
  Hughes, An isogeometric design-through-analysis methodology based on adaptive
  hierarchical refinement of {NURBS}, immersed boundary methods, and {T-spline}
  {CAD} surfaces, {Computer Methods in Applied Mechanics and Engineering} 249
  (2012) 116--150.

\bibitem{hughes2010efficient}
T.~Hughes, A.~Reali, G.~Sangalli, Efficient quadrature for {NURBS}-based
  isogeometric analysis, Computer Methods in Applied Mechanics and Engineering
  199~(5-8) (2010) 301--313.

\bibitem{krishnamurthy2010accurate}
A.~Krishnamurthy, S.~McMains, Accurate moment computation using the {GPU}, in:
  Proceedings of the 14th ACM Symposium on Solid and Physical Modeling, ACM,
  2010, pp. 81--90.

\bibitem{burman2010fictitious}
E.~Burman, P.~Hansbo, Fictitious domain finite element methods using cut
  elements: I. {A} stabilized {Lagrange} multiplier method, {Computer Methods
  in Applied Mechanics and Engineering} 199~(41-44) (2010) 2680--2686.

\bibitem{anderson2018high}
R.~W. Anderson, V.~A. Dobrev, T.~V. Kolev, R.~N. Rieben, V.~Z. Tomov,
  High-order multi-material {ALE} hydrodynamics, SIAM Journal on Scientific
  Computing 40~(1) (2018) B32--B58.

\bibitem{verschaeve2011high}
J.~C. Verschaeve, High order interface reconstruction for the volume of fluid
  method, Computers \& Fluids 46~(1) (2011) 486--492.

\bibitem{dyadechko2005moment}
V.~Dyadechko, M.~Shashkov, Moment-of-fluid interface reconstruction, Los Alamos
  Report LA-UR-05-7571 (2005).

\bibitem{mousavi2010generalized}
S.~Mousavi, H.~Xiao, N.~Sukumar, Generalized {G}aussian quadrature rules on
  arbitrary polygons, International Journal for Numerical Methods in
  Engineering 82~(1) (2010) 99--113.

\bibitem{thiagarajan2014adaptively}
V.~Thiagarajan, V.~Shapiro, Adaptively weighted numerical integration over
  arbitrary domains, Computers \& Mathematics with Applications 67~(9) (2014)
  1682--1702.

\bibitem{scholz2017first}
F.~Scholz, A.~Mantzaflaris, B.~J{\"u}ttler, First order error correction for
  trimmed quadrature in isogeometric analysis, in: Chemnitz Fine Element
  Symposium, Springer, 2017, pp. 297--321.

\bibitem{flusser1999calculation}
J.~Flusser, T.~Suk, On the calculation of image moments, Research Report \#1946
  1946 (1999).

\bibitem{wu2001new}
C.-H. Wu, S.-J. Horng, P.-Z. Lee, A new computation of shape moments via
  quadtree decomposition, Pattern Recognition 34~(7) (2001) 1319--1330.

\bibitem{singer1993general}
M.~H. Singer, A general approach to moment calculation for polygons and line
  segments, Pattern Recognition 26~(7) (1993) 1019--1028.

\bibitem{yang1996fast}
L.~Yang, F.~Albregtsen, Fast and exact computation of {Cartesian} geometric
  moments using discrete {Green}'s theorem, Pattern Recognition 29~(7) (1996)
  1061--1073.

\bibitem{roca2011defining}
X.~Roca, A.~Gargallo-Peir{\'o}, J.~Sarrate, Defining quality measures for
  high-order planar triangles and curved mesh generation, in: Proceedings of
  the 20th International Meshing Roundtable, Springer, 2011, pp. 365--383.

\bibitem{sherwin2002mesh}
S.~Sherwin, J.~Peir{\'o}, Mesh generation in curvilinear domains using
  high-order elements, International Journal for Numerical Methods in
  Engineering 53~(1) (2002) 207--223.

\bibitem{engvall2016isogeometric}
L.~Engvall, J.~A. Evans, Isogeometric triangular {B}ernstein--{B}{\'e}zier
  discretizations: Automatic mesh generation and geometrically exact finite
  element analysis, {Computer Methods in Applied Mechanics and Engineering} 304
  (2016) 378--407.

\bibitem{engvall2017isogeometric}
L.~Engvall, J.~A. Evans, Isogeometric unstructured tetrahedral and
  mixed-element {Bernstein}--{B\'ezier} discretizations, Computer Methods in
  Applied Mechanics and Engineering 319 (2017) 83--123.

\bibitem{engvall2018mesh}
L.~Engvall, J.~A. Evans, Mesh quality metrics for isogeometric
  {B}ernstein--{B}{\'e}zier discretizations, arXiv preprint arXiv:1810.06975
  (2018).

\bibitem{chin2015numerical}
E.~B. Chin, J.~B. Lasserre, N.~Sukumar, Numerical integration of homogeneous
  functions on convex and nonconvex polygons and polyhedra, Computational
  Mechanics 56~(6) (2015) 967--981.

\bibitem{sudhakar2013quadrature}
Y.~Sudhakar, W.~A. Wall, Quadrature schemes for arbitrary convex/concave
  volumes and integration of weak form in enriched partition of unity methods,
  {Computer Methods in Applied Mechanics and Engineering} 258 (2013) 39--54.

\bibitem{saye2015high}
R.~Saye, High-order quadrature methods for implicitly defined surfaces and
  volumes in hyperrectangles, SIAM Journal on Scientific Computing 37~(2)
  (2015) A993--A1019.

\bibitem{olshanskii2016numerical}
M.~A. Olshanskii, D.~Safin, Numerical integration over implicitly defined
  domains for higher order unfitted finite element methods, Lobachevskii
  Journal of Mathematics 37~(5) (2016) 582--596.

\bibitem{sommariva2006meshless}
A.~Sommariva, M.~Vianello, Meshless cubature by {Green}’s formula, Applied
  Mathematics and Computation 183~(2) (2006) 1098--1107.

\bibitem{sommariva2009gauss}
A.~Sommariva, M.~Vianello, {Gauss}--{Green} cubature and moment computation
  over arbitrary geometries, Journal of Computational and Applied Mathematics
  231~(2) (2009) 886--896.

\bibitem{santin2011algebraic}
G.~Santin, A.~Sommariva, M.~Vianello, An algebraic cubature formula on
  curvilinear polygons, Applied Mathematics and Computation 217~(24) (2011)
  10003--10015.

\bibitem{sheynin2003moment}
S.~Sheynin, A.~Tuzikov, Moment computation for objects with spline curve
  boundary, IEEE Transactions on Pattern Analysis and Machine Intelligence
  25~(10) (2003) 1317--1322.

\bibitem{li1993moment}
B.~Li, The moment calculation of polyhedra, Pattern Recognition 26~(8) (1993)
  1229--1233.

\bibitem{jonsson2017cut}
T.~Jonsson, M.~G. Larson, K.~Larsson, Cut finite element methods for elliptic
  problems on multipatch parametric surfaces, {Computer Methods in Applied
  Mechanics and Engineering} 324 (2017) 366--394.

\bibitem{engwirda2014locally}
D.~Engwirda, Locally optimal {D}elaunay-refinement and optimisation-based mesh
  generation (2014).

\bibitem{farin2002curves}
G.~Farin, Curves and surfaces for {CAGD}: A practical guide, Morgan Kaufmann,
  2002.

\bibitem{farouki2012}
R.~T. Farouki, The {B}ernstein polynomial basis: A centennial retrospective,
  Computer Aided Geometric Design 29~(6) (2012) 379--419.
\newblock \href {https://doi.org/10.1016/j.cagd.2012.03.001}
  {\path{doi:10.1016/j.cagd.2012.03.001}}.

\bibitem{davis2007methods}
P.~J. Davis, P.~Rabinowitz, Methods of numerical integration, Courier
  Corporation, 2007.

\bibitem{gautschi1999orthogonal}
W.~Gautschi, Orthogonal polynomials and quadrature, Electron. Trans. Numer.
  Anal 9 (1999) 65--76.

\bibitem{djrbashian1990survey}
M.~Djrbashian, A survey on the theory of orthogonal system and some open
  problems, in: Orthogonal Polynomials, Springer, 1990, pp. 135--146.

\bibitem{gautschi1993Gauss}
W.~Gautschi, Gauss-type quadrature rules for rational functions, in: Numerical
  Integration IV, Springer, 1993, pp. 111--130.

\bibitem{deckers2008rational}
K.~Deckers, J.~Van~Deun, A.~Bultheel, Rational {Gauss}-{Chebyshev} quadrature
  formulas for complex poles outside [-1, 1], Mathematics of Computation
  77~(262) (2008) 967--983.

\bibitem{van2008algorithm}
J.~Van~Deun, K.~Deckers, A.~Bultheel, J.~Weideman, Algorithm 882: Near-best
  fixed pole rational interpolation with applications in spectral methods, ACM
  Transactions on Mathematical Software (TOMS) 35~(2) (2008) 14.

\bibitem{deckers2009computing}
K.~Deckers, J.~Van~Deun, A.~Bultheel, Computing rational {Gauss}--{Chebyshev}
  quadrature formulas with complex poles: The algorithm, Advances in
  Engineering Software 40~(8) (2009) 707--717.

\bibitem{deckers2017algorithm}
K.~Deckers, A.~Mougaida, H.~Belhadjsalah, Algorithm 973: Extended rational
  {F}ej{\'e}r quadrature rules based on {Chebyshev} orthogonal rational
  functions, ACM Transactions on Mathematical Software (TOMS) 43~(4) (2017) 37.

\bibitem{edelman1995polynomial}
A.~Edelman, H.~Murakami, Polynomial roots from companion matrix eigenvalues,
  Mathematics of Computation 64~(210) (1995) 763--776.

\bibitem{MATLAB}
Matlab \texttt{roots} function, the MathWorks, Natick, MA, USA (2019a).

\bibitem{farouki1988algorithms}
R.~T. Farouki, V.~Rajan, Algorithms for polynomials in {Bernstein} form,
  Computer Aided Geometric Design 5~(1) (1988) 1--26.

\bibitem{thiagarajan2018shape}
V.~Thiagarajan, V.~Shapiro, Shape aware quadratures, Journal of Computational
  Physics 374 (2018) 1239--1260.

\bibitem{taber2018moment}
A.~Taber, G.~Kumar, M.~Freytag, V.~Shapiro, A moment-vector approach to
  interoperable analysis, Computer-Aided Design 102 (2018) 139--147.

\bibitem{benzaken2017rapid}
J.~Benzaken, A.~J. Herrema, M.-C. Hsu, J.~A. Evans, A rapid and efficient
  isogeometric design space exploration framework with application to
  structural mechanics, Computer Methods in Applied Mechanics and Engineering
  316 (2017) 1215--1256.

\bibitem{hinz2018spline}
J.~Hinz, M.~M{\"o}ller, C.~Vuik, Spline-based parameterization techniques for
  twin-screw machine geometries, in: IOP Conference Series: Materials Science
  and Engineering, Vol. 425, IOP Publishing, 2018, p. 012030.

\bibitem{kudela2015efficient}
L.~Kudela, N.~Zander, T.~Bog, S.~Kollmannsberger, E.~Rank, Efficient and
  accurate numerical quadrature for immersed boundary methods, Advanced
  Modeling and Simulation in Engineering Sciences 2~(1) (2015) 10.

\bibitem{hermes2018high}
D.~Hermes, P.-O. Persson, High-order solution transfer between curved
  triangular meshes, arXiv preprint arXiv:1810.06806 (2018).

\bibitem{keshavarzzadeh2018numerical}
V.~Keshavarzzadeh, R.~M. Kirby, A.~Narayan, Numerical integration in multiple
  dimensions with designed quadrature, SIAM Journal on Scientific Computing
  40~(4) (2018) A2033--A2061.

\end{thebibliography}

\end{document}